\documentclass[12pt,fullpage,leqno]{article}
\usepackage[english]{babel}
\usepackage{amsmath}
\numberwithin{equation}{section}
\usepackage[margin=1in]{geometry}
\usepackage{amsfonts}
\usepackage{amsmath}
\usepackage{epsfig}
\usepackage{graphicx}
\usepackage{color}
 \usepackage{transparent} 
\usepackage{lmodern} 
\usepackage{mathrsfs}
\usepackage{amsfonts}
\usepackage{amssymb}
\usepackage{mathdots}
\usepackage{amsmath,amssymb,cases,amsthm,color}
\usepackage{float}
\usepackage{subfig}
\usepackage[overload]{empheq}
\usepackage{mathtools}

\usepackage{perpage} 
\MakePerPage{footnote}
\usepackage{fancyhdr} 
\usepackage{indentfirst}
\usepackage{enumerate}
\usepackage{enumitem}
\usepackage[title,titletoc,toc]{appendix}
\usepackage{url}

\usepackage[utf8]{inputenc}
\usepackage[T1]{fontenc}

\usepackage{multirow}
\usepackage{caption}
\usepackage{booktabs}
\usepackage[all]{xy}

\newtheorem{theo}{\bf Theorem}[section]
\newtheorem{lem}[theo]{\bf Lemma}
\newtheorem{pro}[theo]{\bf Proposition}

\newtheorem{defi}[theo]{\bf Definition}
\newtheorem{rem}[theo]{\bf Remark}

\definecolor{cerise}{rgb}{0.87, 0.19, 0.39}

\newenvironment{Proofc}[1]{\smallskip\par\noindent\textsc{#1}\quad}%
  {\hfill$\Box$\bigskip\par}

\def\div{\hbox{div}}

\nonstopmode

\begin{document}
\title{{\bf A solution with free boundary for non-Newtonian fluids with Drucker-Prager plasticity criterion}} 
\author{ Eleftherios Ntovoris\footnote{Université Paris-Est, 
6 et 8 avenue Blaise Pascal 
Cit\'e Descartes - Champs sur Marne 
77455 Marne la Vall\'ee Cedex 2, France.
e-mail: eleftherios.ntovoris@enpc.fr}, M Regis\footnote{
70 rue du Javelot 75013 Paris} } 
\maketitle

\begin{abstract}
We study a free boundary problem which is motivated
by a particular case of the flow of a non-Newtonian fluid, with a pressure depending yield stress given by a Drucker-Prager plasticity criterion. We focus on the steady case and reformulate the equation as a variational problem. The resulting energy has a term with linear growth while we study the problem in an unbounded domain. We derive an Euler-Lagrange equation and prove a comparison principle. We are then able to construct a subsolution and a supersolution which  quantify the natural and expected properties of the solution; in particular we show that the solution has in fact compact support, the boundary of which is the free boundary.

The model describes the flow of a non-Newtonian material on an inclined plane with walls, driven by gravity. We show that there is a critical angle for a non-zero solution to exist. Finally, using the sub/supersolutions we give estimates of the free boundary.
\end{abstract}



\section{Introduction}\label{intro_aval}
\paragraph{Setting of the problem}
We study non-negative solutions $u(y,z)$ of the equation
\begin{equation}\label{eq:intro_E_L}
 \begin{cases} 
\div (\nabla u+|z|q)=-\lambda  & \text{ in } (-1,1)\times (-\infty ,0), \\
 q\in \partial (|\cdot |)(\nabla u), & \text{ }
\end{cases}
\end{equation}
with $u(\pm1,z)=0$, $q=q(y,z)$, $\lambda \geq 0$ and for a function $f:\mathbb{R}^N\rightarrow \mathbb{R}$, $N\in \mathbb{N}$ we define the subdifferential of $f$ at a point $y\in \mathbb{R}^N$ as
\begin{equation}\label{eq:exist_subdiff}
(\partial f)(y):=\{ z\in \mathbb{R}^N:\, f(x)-f(y)\geq z\cdot (x-y)\, \forall x\in \mathbb{R}^N\} .
\end{equation}
The variational formulation of \eqref{eq:intro_E_L} consists in minimizing the functional
\begin{equation}\label{eq:energy_problem}
E_\lambda (u)= \displaystyle\int _{\Omega } \frac{|\nabla u|^2}{2}+|z||\nabla u| -\lambda u ,
\end{equation}
in the space 
\begin{equation}\label{eq:x_space}
\mathcal{X}=\mathcal{X}(\Omega ):=\{u\in W^{1,2} _{0L}(\Omega ),\, z\nabla u\in L^1(\Omega ,\mathbb{R}^2)\} ,
\end{equation}
with
\[
W_{0L} ^{1,2}(\Omega ):= \{ u\in W^{1,2} (\Omega ):\, u(\pm 1,\cdot )=0\} ,
\]
$\Omega =(-1,1)\times (-\infty ,0)$. Note that by Remark \ref{Rem:W^1,1} the functional $E_\lambda$ is well defined in $\mathcal{X}$.
Before we explain the physical interpretation of the mathematical model, we present some of the particularities of the problem. 

Since we study the equation \eqref{eq:intro_E_L} in an unbounded domain, the variational problem \eqref{eq:energy_problem} is no longer trivial because it is not clear if the linear term $-\int_\Omega \lambda u$ is lower semicontinuous or if
 the minimizing sequence obtained by the direct method will have a converging subsequence in $\mathcal{X}$.
Using Lemma \ref{Lem0}, we show that the linear term is lower semicontinuous and the well posedness of the problem is established in Theorem \ref{theo:main_theo_1} \eqref{itm:E_L_eq_1}. Also, despite the fact that the energy $E$ includes a term with linear growth (in the gradient variable), a comparison principle still holds for equation \eqref{eq:intro_E_L}. Using this comparison principle we construct sub/supersolutions and
 show that in fact the solution of \eqref{eq:intro_E_L} is compactly supported. 
 
 For the construction of these barriers we use the “curvature like” equation 
\begin{equation}\label{eq:curv_monn}
-\div (|z|q)=\lambda ,
\end{equation}
which is the first variation of the energy $\int_\Omega |z||\nabla u|-\lambda u$, with $q=\frac{\nabla u}{|\nabla u|}$ when $\nabla u\neq 0$ and $|q|\leq 1$;  then the vector $\frac{\nabla u}{|\nabla u|}$ is the normal to the level sets of $u$. If we suppose that these level sets are given by $-z=\phi (y)$ we are led to study the first variation of the 1-D functional 
\begin{equation}\label{eq:one_dim_cur}
\int_{-1}^1 -\phi (y)\sqrt{1+|\phi '(y)|^2}+\lambda \phi.
\end{equation}
\paragraph{Non-Newtonian fluids}
The model \eqref{eq:intro_E_L} is motivated by the motion of  non-Newtonian fluids. Let $\Omega \subset \mathbb{R}^3$, open and $v:\Omega \rightarrow \mathbb{R}^3$ be the velocity of the fluid, assumed incompressible,
\begin{equation}\label{eq:incompres}
\div\, v=0.
\end{equation}
Let $f:\Omega \rightarrow \mathbb{R}^3$ be the external force, then the relevant equation reads as
\begin{equation}\label{eq:gen_flows}
\div\,\sigma +f=(\nabla\cdot  v)v+\partial _t v
\end{equation}
where $\sigma $ is the stress tensor and using the usual summation convention we write $(\nabla\cdot  v)v=(v_j\partial _jv_i)_{1\leq i\leq 3}$.
Let $\sigma _{\mathrm{dev}}$ be the stress deviator defined by $\mathrm{tr}(\sigma _{\mathrm{dev}})=0$ and
\begin{equation}\label{eq:dev_stres}
\sigma _{\mathrm{dev}} :=\sigma +pI,
\end{equation}
where $p$ is the pressure and $I$ is the unit matrix.\par
\begin{figure}[htbp]
     \centering
     \includegraphics[scale=1]{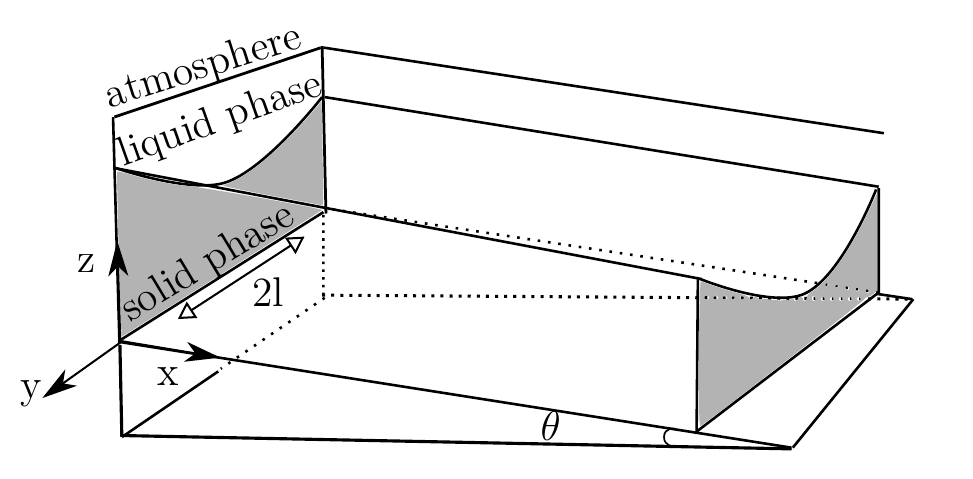}
     \caption{Steady flow on an inclined plane.}
    \label{fig:intro_aval}
\end{figure}

 We are interested in the flow of rigid visco-plastic fluids, which unlike Newtonian fluids can sustain shear stress. The stress tensor in this case is characterized by a flow/no flow condition, namely when the stress tensor belongs to a certain convex set the fluid behaves like a rigid body, whereas outside this set the material flows like a regular Newtonian fluid. For a matrix $B=(b_{ij})_{1\leq i,j\leq 3}$ we denote the norm $||B||=\sqrt{\frac{1}{2}\displaystyle\sum_{i,j=1}^3 b_{ij} ^2}$.
  Following \cite{Ionescu20151} and \cite{Cazacu20061640} we define the stress deviator as
 \begin{equation}\label{eq:constit_dr_pr}
 \begin{cases} 
\sigma _{\mathrm{dev}}=2\nu D(v)+k(p)\frac{D(v)}{||D(v)||}  & \text{ if } D(v)\neq 0, \\
 ||\sigma _{\mathrm{dev}}||\leq k(p) & \text{ if } D(v)=0
\end{cases}
 \end{equation}
 where we assume that the viscosity $\nu >0$ is constant and $k(p)$ is the pressure-dependent yield stress and $D(v)=(\nabla v+(\nabla v)^T)/2$. The above constituent law is a result of a superposition of the viscous contribution $2\nu D(v)$ and a contribution related to plasticity effects $k(p)\frac{D(v)}{||D(v)||}$, which is independent of the norm of the strain rate $||D(v)||$. For constant yield limit $k(p)$ we retrieve the regular Bingham model, which is a generalized Newtonian problem, i.e. the constituent law in this case is described by a dissipative potential, see \cite{MR0521262}, \cite[Chapter 3]{MR1810507} and references therein. In this paper we will assume the Drucker-Prager plasticity criterion
 \begin{equation}\label{eq:druck_prag}
 k(p)=\mu _sp,
 \end{equation}
where $\mu _s=\tan \delta _s$, with $\delta _s$ the internal friction (static) angle. The existence of a dissipative potential in the case of Bingham flows allows for a variational formulation and in tern the well-posedeness of the problem; for quasi-static Bingham flows see for example \cite{MR1810507}. The case of a Drucker-Prager criterion, however, falls in a wider class of constituent laws called “$\mu(I)-$rheology” which are known to be ill-posed, see \cite{FLM:9921740} and \cite{SCHAEFFER198719}. The strong geophysical interest in the model \eqref{eq:druck_prag} supports however our study. A main result of the present work is that for one-directional steady flows the model is well-posed.
\paragraph{Flow in one direction} We study the well-posedeness and certain quantitative properties of quasi-static solutions  of \eqref{eq:incompres}-\eqref{eq:gen_flows}, \eqref{eq:constit_dr_pr}-\eqref{eq:druck_prag}, for a material which flows on an inclined plane with sidewalls. We assume that the inclination angle is constant $\theta$ and the material moves only in the direction $x$ under the effect of gravity, see Figure \ref{fig:intro_aval}.
In what follows, we will assume that the velocity field is of the form $v(x,y,z)=(u(y,z),0,0)$ for $(x,y,z)\in\{(x,y,z):\, 0\leq x,\, -l\leq y\leq l,\, b\leq z\leq h(y)\}$ where $h(y)$ is the interface separating the fluid and the air and $z=b$ is the surface of the inclined plane, the width of which is equal to $2l$. Although the well posedness of similar problems have been studied in more generality in a bounded domain, as it will become clear later, in order to study the interface between the solid and the liquid phase, as we increase the inclination angle, we will need to take $b=-\infty$. By the form of $v$, the incompressibility condition \eqref{eq:incompres} is trivially satisfied and
equation \eqref{eq:gen_flows} with $\partial _t v=0$ becomes
\begin{equation}\label{eq:quasi_static_1}
\div \sigma _{\mathrm{dev}}=-f+\nabla p,
\end{equation}
with $\sigma _{\mathrm{dev}}$ given by \eqref{eq:constit_dr_pr} and 
\begin{equation}\label{eq:gravity}
f=(g_0\sin\theta ,0,-g_0\cos\theta ),
\end{equation}
with $g_0$ the gravitational constant.
We also assume that $p=p(y,z)$. We calculate
\begin{equation}\label{eq:rate_strain_vel}
D(v)=\frac{1}{2}\left( \begin{array}{ccc}
0 & \partial _y u & \partial _z u \\
\partial _y u & 0 & 0 \\
\partial _z u & 0 & 0 \end{array} \right)
\end{equation}
and $||Du||=\frac{1}{2}|\nabla u|$, with $\nabla u=(\partial _yu,\partial _z u)$. If we substitute \eqref{eq:rate_strain_vel} in \eqref{eq:constit_dr_pr}, equations \eqref{eq:quasi_static_1} become, for $D(v)\neq 0$ or equivalently $\nabla u\neq 0$,
\begin{alignat}{3}[left=\empheqlbrace]
  \nu \div (\nabla u)+\mu _s\div \left( p\frac{\nabla u}{|\nabla u|}\right) &=  -g_0\sin \theta \label{eq:first_quasi_static}\\
  0 &= \partial _y p \label{eq:der_pres_y} \\
  0&= g_0\cos \theta +\partial _z p\label{eq:pressure}
\end{alignat}
Where the divergence is taken for the coordinates $(y,z)$. If we integrate equation \eqref{eq:pressure} from $z$ to $h(y)$ we get $p(y,z)=(h(y)-z)g_0\cos \theta ,$
but because of equation \eqref{eq:der_pres_y} and because $\theta \in \left[ 0,\frac{\pi}{2}\right)$ we have $h(y)=h\equiv \text{ constant}$. For simplicity we take $h=0$; then the pressure is given by
\begin{equation}\label{eq:pressure_funct}
p(y,z)=|z|g_0\cos \theta .
\end{equation}
We are lead to study the following equation
 \begin{equation}\label{eq:single_eq_prob}
 \begin{cases} 
\nu \,\div (\nabla u)+\mu _s g_0\cos\theta \,\div (|z|q)=-g_0\sin \theta  & \text{ in } (-l,l)\times (b,0), \\
 q\in \partial (|\cdot |)(\nabla u) & \text{ }
\end{cases}
 \end{equation}
where $\partial (|\cdot |)$ is the subdifferential of the absolute value. If $(u,q)$ is such that \eqref{eq:single_eq_prob} holds, with $q=q(y,z)=(q_1(y,z),q_2(y,z))$, then $|q|\leq 1$ and $q=\frac{\nabla u}{|\nabla u|}$ for $\nabla u \neq 0$ and therefore the stress deviator defined by
\begin{equation}\label{eq:proj_stress}
\sigma _{\mathrm{dev}}:=\nu \left( \begin{array}{ccc}
0 & \partial _y u & \partial _z u \\
\partial _y u & 0 & 0 \\
\partial _z u & 0 & 0 \end{array} \right)
+\mu _s |z|g_0\cos \theta \left( \begin{array}{ccc}
0 & q_1 & q_2 \\
q_1 & 0 & 0 \\
q_2 & 0 & 0 \end{array} \right)
\end{equation}
is of the form \eqref{eq:constit_dr_pr} with $v(x,y,z)=(u(y,z),0,0)$ and solves equations \eqref{eq:quasi_static_1} with $f$ given by \eqref{eq:gravity} and $p$ by \eqref{eq:pressure_funct}.
\paragraph{Boundary conditions}
On the surface of the material $z=0$ we assume a \emph{no stress condition}, i.e. $\sigma \cdot (0,0,1)=0$; since the pressure is zero on the surface near the atmosphere, this condition becomes $\sigma _{\mathrm{dev}}\cdot (0,0,1)=0$. Here we assume that the stress deviator is given by \eqref{eq:proj_stress}. Then the stress free condition becomes (since $z=0$)
\begin{equation}\label{eq:stress_free_cond}
\partial _zu (y,0)=0. 
\end{equation}
On the lateral boundary $y=\pm 1$ we assume the Dirichlet conditions $u=0$ (no slip), while at the bottom $z=b$, where the material is in contact with the inclined plane, a natural assumption is the friction condition
\[ 
 \begin{cases} 
\sigma n-(\sigma n\cdot n)n =\mu _Cv  & \text{  }  \\
v\cdot n =0 & \text{ }
\end{cases}
\]
where $v,\sigma ,n,\mu _C$ are the velocity, stress, normal to the plane and a friction coefficient respectively. In our case the friction condition reads as follows
\begin{equation}\label{eq:friction_cond}
\nu \partial _zu +\mu _s |b|g_0(\cos\theta )q_2 =\mu _Cu.
\end{equation}
\paragraph{Variational formulation}
 The variational formulation of equation \eqref{eq:single_eq_prob} with boundary conditions \eqref{eq:stress_free_cond}, \eqref{eq:friction_cond} and the homogeneous Dirichlet conditions on the lateral boundary constitutes in minimizing the energy
\begin{equation}\label{eq:phys_energy}
\int_{(-l,l)\times (b,0) } \nu \frac{|\nabla u|^2}{2}+\mu _s |z|g_0\cos\theta |\nabla u|-(g_0\sin\theta ) u +\mu _C\int_{\{z=b\}}\frac{|u|^2}{2}
\end{equation}
with zero lateral boundary conditions, i.e. $u(\pm 1,\cdot )=0$. Since the energy \eqref{eq:phys_energy} is convex and the domain is bounded we can easily get a non-negative minimizer via the direct method.\par
 We are interested in the properties of the minimizer as we increase the inclination angle $\theta$.
We call solid and liquid phases the sets $\{(y,z): u(y,z)=0\}$ and $\{(y,z): u(y,z)>0\}$ respectively (often abbreviated as $\{u=0\}$, $\{u>0\}$ resp.), while their common boundary we call a yield curve. We note that usually in the literature the yield curve is defined, for our setting, as the set $\partial \{\nabla u\neq 0\}$, but approximating this set would require different methods and more regularity of the solution.

For $|b|$ small we expect that for a sufficiently large angle $\theta$ all of the material will move due to the gravity, namely there is no solid phase, whereas, if $|b|$ is large enough, even if the inclination is large we expect that there will be a solid phase. In order to study the behaviour and shape of the liquid/solid phases as we increase the inclination angle, we fix $b=-\infty$. However, there is still one more free boundary remaining, the yield curve, i.e. the curve that separates the solid from the liquid phase. Since we study \eqref{eq:phys_energy} in an unbounded domain we drop the friction condition.
Let $\tilde{u}$ be a solution of \eqref{eq:single_eq_prob}-\eqref{eq:stress_free_cond} with $b=-\infty$, in order to simplify further the equation \eqref{eq:single_eq_prob} we set 
\begin{equation}\label{eq:change_mass}
u(y,z)=\frac{\nu}{\mu _sg_0\cos\theta}\frac{\tilde{u}\left( ly,lz\right)}{l^2} , \quad (y,z)\in (-1,1)\times (-\infty ,0) 
\end{equation}
 we also define
\begin{equation}\label{eq:abs_rel_lagr_mult}
\lambda := \frac{\tan \theta}{\mu _s},
\end{equation}
then $ \partial |\cdot |(\nabla u(y,z))=\partial |\cdot |(\nabla \tilde{u}(ly,lz))$ and therefore, $u$ given by \eqref{eq:change_mass} solves the equation \eqref{eq:intro_E_L}
if and only if $\tilde{u}$ solves \eqref{eq:single_eq_prob} 

 As we will see in Theorem \ref{Pro:main_prop}, the minimizer has compact support, therefore, it trivially satisfies the friction condition \eqref{eq:friction_cond} on the solid phase as long as the level of the plane is taken far enough from the support of the minimizer.\par
\begin{figure}[htbp]
     \centering
     \includegraphics[scale=0.8]{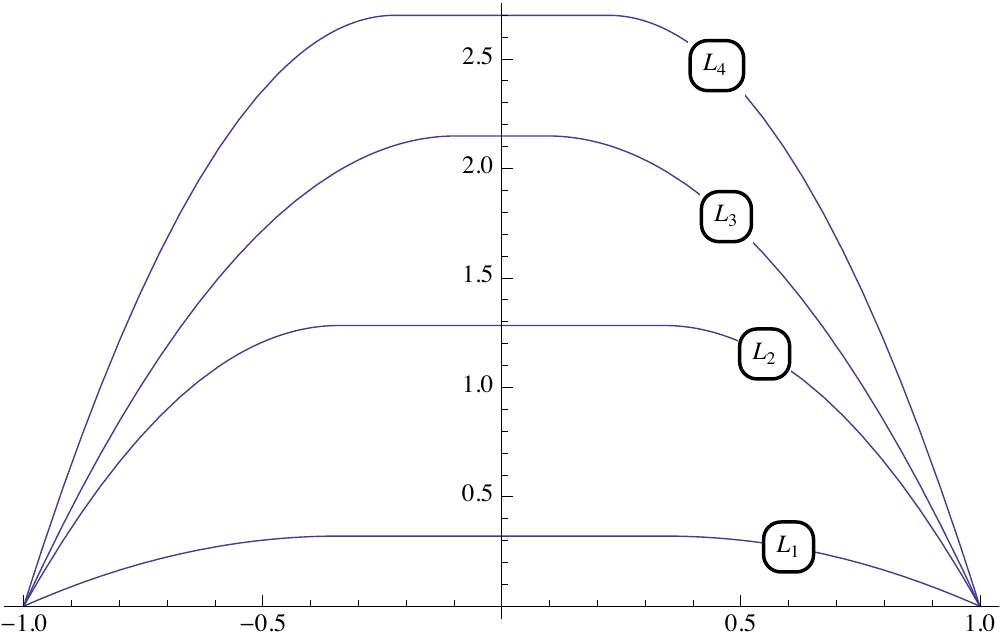}
     \caption{$L_i,\, i=1,2,3,4$ are the graph of $w$ for $(m ,A)=(0.5,0.5),(2,2),(3,0.5),(4,2)$ respectively.}
    \label{fig:1D_minimiser}
\end{figure}

We also show that the critical angle for an non-zero minimizer to exist is $\arctan\mu _s$, namely for $\theta >\arctan\mu _s$ there exists a non-zero solution with a yield curve while for $0\leq \theta \leq \arctan\mu _s$ the solution is zero. This angle
is known in the literature by experimental study, see for example \cite{1742-5468-2006-07-P07020}. The time dependent, one dimensional analogue of our case is studied in \cite{Bouchut}; the authors prove that for $\theta >\arctan\mu _s$ there is no solution with solid phase while in our case the solution always has a solid phase. The difference of course lies in our two dimensional setting of the problem in which the existence of the walls where the velocity vanishes is crucial, not just for the physical relevance of the problem. Indeed since we study minimizers of \eqref{eq:energy_problem} in an unbounded domain we will often need to apply Poincar\'e's inequality, for this reason we need that the projection of the domain in one of the coordinate axes is bounded. In \cite{mang} the authors also prove that for $\theta \leq\arctan\mu _s$ the flowing material stops moving in finite time.
\paragraph{Review of the literature}
For an extensive review of non-Newtonian fluids see \cite{MR0521262}, also \cite{MR1810507} and references therein and \cite{MR1409366} for evolutionary problems. The flow of a viscoplastic material with “$\mu (I)-$rheology” is relatively new in the literature, see for example \cite{Ionescu20151}. The inviscid case, i.e. for $\nu =0$ is similar to another scalar model with applications in image processing, the total variation flow, see for example \cite{MR2033382} and \cite{MR1929886}. Although the total variation bears more similarities with the Bingham case, many of the tools used to study our problem are similar. In fact the total variation is more difficult to study because of the lack of the quadratic term in the energy which leads to lack of regularity of the solution. For the inviscid case our energy \eqref{eq:energy_problem} falls into a wider class, the “total variation functionals” see \cite[Hypothesis 4.1]{gen_consti}. We refer to \cite{Roche} for simulations of a regularized Drucker-Prager model with application to granular collapse. Concerning the case of the inclined plane see \cite{Jop}  and \cite{1742-5468-2006-07-P07020}.
\paragraph{Organization of the paper}
In Section \ref{Sec:main_res} we state our main results, Theorems \ref{theo:main_theo_1} and \ref{Pro:main_prop}. In Subsection \ref{SubSec:1d} we study the 1-dimensional analogue of \eqref{eq:energy_problem} which we use in Lemma \ref{Lem0}; this Lemma is the crucial step in order to prove that the linear term $-\lambda \int_\Omega u$ is lower semicontinuous. In Subsection \ref{Subsec:Regul} we study an approximate problem of the minimizer of \eqref{eq:energy_problem} which helps us to prove certain regularity properties of the solution; we also note that since the minimizer is studied in the half stripe $\Omega$ the regularity holds up to the interface seperating the solid from the liquid phase (the support of the minimizer). Using the approximate minimizer we can also calculate the first variation of \eqref{eq:energy_problem}. Finally, in Lemma \ref{Lem:explicit_minim} we construct a solution of \eqref{eq:curv_monn} which we use together with the comparison principle from Subsection \ref{SecComparison}, in Subsections \ref{SecSubsolution} and \ref{SecSupersolution} in order to construct a subsolution and supersolution respectively. The Figures \ref{fig:1D_minimiser}-\ref{fig:increase_min} as well as the simulations in Table \ref{tab:test1} have been made with Mathematica.
\section{Main results}\label{Sec:main_res}


We begin with a technical remark.
 \begin{rem}\label{Rem:W^1,1}
We have $\mathcal{X}(\Omega )\subset W^{1,1}(\Omega )$, which justifies the choice of the space $\mathcal{X}$ as natural functional space for the functional \eqref{eq:energy_problem}. Indeed,
\[
\int_{\{|z|\geq 1\}\cap\Omega } |\nabla u|\leq \int_{\{|z|\geq 1\}\cap\Omega }|z||\nabla u|<\infty ,
\]
from which get that $u\in L^1(\Omega )$ by Poincar\'e's inequality, see \cite[Theorem 12.17]{MR2527916}; note also that in our case the proof of Poincar\'e's inequality requires only that elements of the space $W^{1,2} (\Omega )$ are zero on the lateral boundary of $\Omega$ (i.e. on $\{\pm 1\}\times (0,\infty )$). In fact, since the width of the walls is $2$ we have $\int_\Omega |u|^p\leq \frac{2^p}{p}\int_\Omega |\nabla u|^p$ for $p=1,2$.
\end{rem}


Let
\begin{equation}\label{eq:lambda_def}
\Lambda :=\lbrace q:\,  q\in  L_{loc} ^2(\Omega ,\mathbb{R}^2),\, |q|\leq 1 \, \text{ a.e.}\rbrace .
\end{equation}
Let $\hat{\Omega}=(-1,1)\times \mathbb{R}$, $u\in W^{1,2} _{0L}(\Omega )$, we denote by $\hat{u} \in W^{1,2} _{0L} (\hat{\Omega})$ the reflection of $u$ with respect the $z=0$ axes, i.e.
\begin{equation}\label{eq:reflex}
\hat{u} (y,z):=
\begin{cases} 
u(y,z)& \text{ if } (y,z)\in \Omega ,\\
u(y,-z)   &  \text{ if }y,z)\in \hat{\Omega}\setminus \Omega .
\end{cases} 
\end{equation}
Throughout the paper we will denote the space $\mathcal{X}(\Omega )$ simply by $\mathcal{X}$. Only in Lemma \ref{Lem1} we will use the explicit notation, this time for the space $\mathcal{X}(\hat{\Omega})$.
The weak formulation of \eqref{eq:intro_E_L} is 
 \begin{equation}\label{eq:sub_differential}
 \begin{cases} 
 \displaystyle\int_\Omega \nu \nabla u\cdot \nabla \varphi +|z| q\cdot \nabla \varphi =\lambda\int_\Omega  \varphi   & \text{for all } \varphi \in\mathcal{X} \\
q\cdot \nabla u=|\nabla u| & \text{a.e.}
\end{cases}
 \end{equation}
 for some $\lambda \geq 0$, $q\in \Lambda$.
We can now state our first main Theorem.
\begin{theo}\label{theo:main_theo_1}{\bf(Existence and uniqueness of minimizers of \eqref{eq:energy_problem})}\\
Let $\lambda \geq 0$, $E_\lambda$ be given by \eqref{eq:energy_problem}, then the following hold
\begin{enumerate}[label=(\roman*) ,ref=\roman*]
\item \label{itm:E_L_eq_1} there exists a unique $0\leq u_{\lambda } \in \mathcal{X}  $ such that
\begin{equation}\label{eq:min_prob_1}
E_\lambda (u_\lambda )=\inf_{v\in \mathcal{X} }E_\lambda (v) ,
\end{equation}
moreover, $u_\lambda \equiv 0$ if $\lambda \in [0,1]$ and $u_\lambda \not\equiv 0$ if $\lambda \in (1,+\infty )$,
\item \label{itm:E_L_eq_2}there exists $q\in \Lambda $ such that $(u_\lambda ,q )$ solves \eqref{eq:sub_differential},
\item \label{itm:E_L_eq_3}$u_\lambda \in C^{0,\alpha } _{\mathrm{loc}}(\Omega )$ for all $\alpha \in (0,1)$, in fact $\hat{u} _\lambda\in W^{2,2} _{\mathrm{loc}}(\hat{\Omega})$ and $\partial _z u_\lambda (y,0)=0$ for $y\in (-1,1)$,
\item \label{itm:E_L_eq_4}if $\lambda >1$, the pair $(u_\lambda ,q )$ obtained in \eqref{itm:E_L_eq_2} is unique in the sense that if $(\bar{u}_\lambda ,\bar{q})\in \mathcal{X }\times \Lambda$ is another pair satisfying \eqref{eq:sub_differential} then
\[
	u=\bar{u} \text{ in }\Omega , \,\,\text{ and } q=\bar{q} ,\text{ a.e. in } \{\nabla u\neq 0\} .
	\]
\end{enumerate}
\end{theo}
We set 
\[
I_{m}:=\inf_{v\in \mathcal{X} }E_\lambda (v).
\]
Note that by the continuity of the non-negative function $u_\lambda $ in Theorem \ref{theo:main_theo_1} we can define the \emph{yield curve} as the common boundary $\partial \{u_\lambda >0\}=\partial \{u_\lambda =0\}$. Moreover, the critical value $\lambda =1$ in the previous Theorem is also a critical value of the physical solution by \eqref{eq:change_mass}, \eqref{eq:abs_rel_lagr_mult} and it does not depend on the viscosity constant $\nu$ or the width of the walls.

We will give some notations in order to present our second result, the motivation for this notation will become clear in the proofs of the relevant Propositions. Let $\lambda >1$ for $Z\in [\frac{1}{\lambda },\frac{1}{\lambda -1}]$ we define
\begin{equation}\label{eq:inverce_funct}
f_\lambda (Z):=
 \frac{1}{(\lambda ^2-1)^{3/2}}\left\lbrace \mathrm{Arcsin}\left[(\lambda ^2-1)Z -\lambda \right]-\lambda \sqrt{1-\left( (\lambda ^2-1)Z-\lambda \right) ^2} \right\rbrace .
\end{equation}
As we will see in the proof of Lemma \ref{Lem:explicit_minim}, the function $f_\lambda$ is strictly increasing in the interval $[\frac{1}{\lambda},\frac{1}{\lambda -1}]$, i.e. $f_\lambda (Z_1)<f_\lambda (Z_2)$ for $Z_1<Z_2$, with $Z_1,Z_2\in[\frac{1}{\lambda},\frac{1}{\lambda -1}]$; we can therefore define the following function
\begin{equation}\label{eq:expl_implicite_solution}
\phi _{K(\lambda )}(y):=
K(\lambda ) f_\lambda ^{-1} \left( f_\lambda \left( \frac{1}{\lambda -1}\right)+\frac{|y|}{K(\lambda )}\right) \quad \text{ } y\in [-1,1],
\end{equation}
where
\begin{equation}\label{eq:K_max}
K(\lambda ):=\frac{1}{f_\lambda \left( \frac{1}{\lambda}\right) -f_\lambda \left( \frac{1}{\lambda -1}\right)} .
\end{equation}
Note that by the monotonicity of $f_\lambda$ it is $K(\lambda )<0$.
We also define the half cone 
\begin{equation}\label{eq:cone}
\mathcal{C}_\lambda := \{ (y,z)\in \mathbb{R}^2:\,0<|y|< z\frac{\lambda}{K(\lambda )}\} 
\end{equation}
 and 
 \begin{equation}\label{eq:Epi}
\mathrm{Epi} ^\smallsmile (\lambda ):=\{(y,z)\in \Omega: z>\phi _{K(\lambda )}(y)\} .
 \end{equation}
 In Lemma \ref{Lem:explicit_minim} we show that the sets in \eqref{eq:Epi} are increasing in $\lambda$ in the sense that $\mathrm{Epi} ^\smallsmile (\lambda )\subsetneq \mathrm{Epi} ^\smallsmile (\bar{\lambda})$ for $\bar{\lambda}> \lambda$, see Figure \ref{fig:lambdas}.
For $\lambda _1>\lambda $ we set
\begin{equation}
\vartheta _{\lambda ,\lambda _1}:=\frac{\lambda _1-\lambda }{2 \left(1+\left( \frac{\lambda _1}{K(\lambda _1)}\right) ^2\right)} ,\label{eq:theta}
\end{equation}
\begin{equation}
b(\lambda ,\lambda _1):=1+\sqrt{\frac{\lambda _1}{2\vartheta _{\lambda ,\lambda _1}}},\label{eq:b}
\end{equation}
\begin{equation}
\Pi (\lambda ,\lambda _1):=\frac{-K(\lambda _1)}{\lambda _1-1}b(\lambda ,\lambda _1)+\frac{K(\lambda )}{\lambda -1}.\label{eq:Pi}
\end{equation}
and
\begin{equation}\label{eq:Sup}
\mathrm{Epi} _\smallsmile (\lambda _1):=\lbrace (y,z)\in \Omega :\, z>b(\lambda ,\lambda _1)\phi _{K(\lambda _1)}\left(\frac{y}{b(\lambda ,\lambda _1)}\right) \rbrace .
\end{equation}
In Lemma \ref{Lem:explicit_minim} we see that $\displaystyle\min_{|y|\leq 1}\phi _{K(\lambda )}(y)=\phi _{K(\lambda )}(0)=\frac{K(\lambda )}{\lambda -1}$ for all $\lambda >1$, and therefore, the function $\Pi$ in \eqref{eq:Pi} is the distance of the projections on the $z-$axes of the epigraphs $\mathrm{Epi} ^\smallsmile (\lambda )$ and $\Omega\setminus \mathrm{Epi} _\smallsmile (\lambda _1)$.
Using \eqref{eq:K_max} we calculate
\[
\frac{K(\lambda _1)}{\lambda _1-1}=\frac{2(\lambda _1+1)\sqrt{\lambda _1 ^2-1}}{2\sqrt{\lambda _1 ^2-1}+\pi +2\mathrm{Arcsin}\left( \frac{1}{\lambda _1}\right)} ,
\]
then $\displaystyle\lim_{\lambda _1\rightarrow +\infty}\frac{K(\lambda _1)}{\lambda _1-1}=+\infty$, and similarly one can see that $\displaystyle\lim_{\lambda _1\rightarrow +\infty}\frac{K(\lambda _1)}{\lambda _1}=+\infty$; if we combine the above two limits, one can check that for all $\lambda >1$, 
\begin{equation}\label{eq:limit_1dire}
\displaystyle\lim_{\lambda _1\rightarrow+\infty} \Pi (\lambda ,\lambda _1)=+\infty .
\end{equation}
 We also have
\begin{equation}\label{eq:limit_2dire}
\displaystyle\lim_{\lambda _1\rightarrow \lambda}\Pi (\lambda ,\lambda _1)=+\infty .
\end{equation} 
If we combine \eqref{eq:limit_1dire}, \eqref{eq:limit_2dire} and the fact that $\Pi$ is continuous and we get that for every fixed $\lambda >1$ the function $\Pi (\lambda ,\cdot )$ attains a minimum for some $\lambda ^\star _1>\lambda$.
In fact numerical simulations (see Figure \ref{fig:test2}) suggest the function $\Pi (\lambda ,\cdot )$ attains the minimum at a unique $\lambda ^\star _1>\lambda $, but the analytical calculations are too complicated to check.\par
In the following Theorem we gather the main properties of the solution obtained in Theorem \ref{theo:main_theo_1}.
\begin{theo}{\bf (Main properties)\\}\label{Pro:main_prop}
Let $\lambda >1$, $\mathcal{X}$ as in \eqref{eq:x_space}, $(u _\lambda ,q)\in \mathcal{X}\times \Lambda$ be a solution of \eqref{eq:sub_differential}. Also let $\lambda _1>\lambda$, $\mathrm{Epi} _\smallsmile (\lambda _1)$ be as in \eqref{eq:Sup}, then then function $u_\lambda $ has compact support and it's support can be estimated as follows
	\begin{equation}\label{eq:free_bd_estimate}
	\mathrm{Epi} ^\smallsmile (\lambda )\subset \mathrm{supp}\, u_\lambda \subset \mathrm{Epi} _\smallsmile (\lambda _1) .
	\end{equation}
	Moreover, we can optimize estimate \eqref{eq:free_bd_estimate} by choosing $\lambda _1=\lambda ^\star _1$.
\end{theo}
\begin{rem}{\bf (Consequences of Theorem \ref{Pro:main_prop})}
\begin{enumerate}
\item In Lemma \ref{Lem:explicit_minim} we show that the function $\phi _{K(\lambda )}$ has a strictly negative maximum, therefore estimate \eqref{eq:free_bd_estimate} implies that the yield curve $\partial \{ u_\lambda >0\}$ never reaches the surface of the atmosphere $\{z=0\}$.
\item Notice that the sets $\mathrm{Epi} ^\smallsmile (\lambda )$ and $\mathrm{Epi} _\smallsmile (\lambda _1)$ can also estimate the support of the physical solution, by \eqref{eq:change_mass} and they are independent of the viscosity $\nu$.
\end{enumerate}
\end{rem}
\section{Existence/Uniqueness}

\subsection{$1D$-problem}\label{SubSec:1d}
Let $A>0$ and for $w\in W^{1,2} _0(-1,1) $ we consider the energy
\begin{equation}\label{ene-1D}
\epsilon _{A}(w)=\int_{-1}^{1}\left(\frac{ |w'(y)|^{2}}{2}+A|w'(y)|\right)dy .
\end{equation}

Using the direct method of calculus of variations it is not difficult to show the following Proposition.
\begin{pro}\label{min-1D}{\bf(Minimizer of $1D$-problem)}\\
Let $A>0$ and $m>0.$ Then there exists a unique function $w$ solving 
\[
\epsilon _A(w)=\inf_{\substack{\overline{w}\in W^{1,2} _0(-1,1)\\ \int_{-1}^{1}\overline{w}=m}} \epsilon _A(\overline{w}).
\]
\end{pro}
We set 
\begin{equation}\label{min}
I^{A}_{m}:=\inf_{\substack{\overline{w}\in W^{1,2} _0(-1,1)\\ \int_{-1}^{1}\overline{w}=m}} \epsilon _A(\overline{w}).
\end{equation}
The uniqueness of the minimizer of $\epsilon _A$ in the above Proposition follows by the strict convexity of the functional or by using similar arguments as in the proof of Step 1 of Theorem \ref{theo:main_theo_1} \eqref{itm:E_L_eq_1}.\par
We define the set theoretic sign function as
\[ \mathrm{sign}(r):=
\begin{cases} 
\left\lbrace \frac{r}{|r|}  \right\rbrace & \text{ if }  r\neq 0 ,\\
 (-1,1)  &  \text{ else.}
\end{cases} 
 \]
\begin{pro}\label{chara-1D-sol}{\bf(Characterization of the $1D$ minimizer)}\\
Let $A,\, m >0$. If $\lambda _A=\lambda _{A,m}$ is the non-negative root of
\begin{equation}\label{eq:1D_multiplier}
2\lambda _A ^3-3\lambda _A ^2 (A+m )+A^3=0,
\end{equation}
with $\lambda _A> A+m$,
\begin{equation}\label{eq:1D_width}
a=\frac{A}{\lambda _A}< 1,
\end{equation}
\begin{equation}\label{form-w}
w(y)=\begin{cases} 
A\left( -\frac{y^2}{2a}+|y|+\frac{1}{2a}-1\right)& a<|y|< 1,\\\\
A\frac{(a-1)^2}{2a} &  |y|\leq a.
\end{cases} 
\end{equation}
 \[ 
q(y)=\begin{cases} 
-\frac{y}{|y|}& a<|y|< 1,\\\\
 -\frac{y}{a}  &  |y|\leq a,
\end{cases} 
\]
then $(w,q,\lambda _A)$ solves the equation 
\begin{equation}\label{eq:1DEuler_lagr}
-w''(y)-A(q(y))'=\lambda _A,\quad \text{for a.e. } y\in (-1,1),
\end{equation}
 and $\int_{-1}^1w=m $. In particular $w$ is the unique minimizer of \eqref{ene-1D} corresponding to the volume constraint $m$.
\end{pro}
\begin{Proofc}{\bf Proof of Theorem \ref{chara-1D-sol}}\\
\textbf{Step 1. The trinomial \eqref{eq:1D_multiplier}}\par
First we will show that the trinomial \eqref{eq:1D_multiplier} has a unique non-negative root $\lambda _A\geq A+m$.  A simple calculation shows that the trinomial \eqref{eq:1D_multiplier} is increasing in the interval $[A+m,+\infty )$ with values $A^3\left( 1-\left( 1+\frac{m }{ A}\right) ^3\right)\leq 0$ and $+\infty$ for $\lambda _A=A+m$ and $\lambda _A=+\infty$ respectively, hence there exist a root for $\lambda \geq A+m$.\\
\textbf{Step 2. The equation \eqref{eq:1DEuler_lagr}}\par
For a.e. $y\in (-1,1)$ we have
\begin{equation}\label{eq:1D_deriv}
w''(y)=\begin{cases} 
-\frac{A}{a} & a<|y|< 1,\\
0 &  |y|< a
\end{cases} 
\end{equation}
and
\begin{equation}\label{eq:1D_deriv_subdiff}
q'(y)=\begin{cases} 
0 & a<|y|< 1,\\
-\frac{1}{a} &  |y|< a.
\end{cases} 
\end{equation}
Using \eqref{eq:1D_deriv}, \eqref{eq:1D_deriv_subdiff} and \eqref{eq:1D_width} we deduce that $(w,q,\lambda _A)$ solves \eqref{eq:1DEuler_lagr}.\\
\textbf{Step 3. Volume constraint}\par
 It remains to show that $\int_{-1}^1 w=m $. It is
\begin{align*}
\int_{-1}^1 w &=2\left( \frac{A}{2}(1-a)^2 +\int_a^1 w\right) \\
&= 2\left( \frac{A}{2}(1-a)^2 -\frac{A}{2a}\int_a^1 (y-a)^2-(1-a)^2\, dy \right)\\
&= \frac{A(1-a)^2}{a}-\frac{A}{a}\int_0^{1-a}y^2\,dy\\
&=\frac{A}{3a}(1-a)^2(2+a),
\end{align*}
using equation \eqref{eq:1D_width} and \eqref{eq:1D_multiplier} we get
\begin{equation}\label{mu-A}
\int_{-1}^1w=\frac{A}{3a}(1-a)^2(2+a)=\frac{(\lambda _A -A )^2(2\lambda _A +A)}{3\lambda _A ^2}=m.
\end{equation}
\textbf{Step 3. Minimizer}\par
It remains to show that $w$ is the minimizer of $\epsilon _A$ in $W^{1,2} _0 (-1,1)$ which corresponds to the constraint $m$. First we notice that $q(x)\in \mathrm{sign}(w'(x))=\partial (|\cdot |)(w'(x))$ for $x\in (-1,1)$ and the subdifferential is given by \eqref{eq:exist_subdiff}.
Let $v\in W^{1,2} _0(-1,1)$ with $\int_{-1}^1 v=m$, it is
\begin{align*}
\epsilon _A(v)-\epsilon _A(w)& \geq \int_{-1}^1 w'(v-w)'+Aq(v-w)'\\
&= -\int_{-1}^1 (w''+Aq')(v-w)=\lambda _A\int_{-1}^1(v-w)=0.
\end{align*}
\end{Proofc}
\subsection{A variational problem}\label{SecExistence}
The lower semicontinuity of the term $-\int_\Omega \lambda u$ in \eqref{eq:energy_problem} under the weak topology of $W^{1,2}$ is not trivial since the integral is not evaluated in a bounded domain.
The following Lemma shows that the $L^1$-tails of a sequence of functions will converge to zero if the respective values of the functional $E_{\lambda}$ are uniformly bounded.
\begin{lem}{\bf(Compensation of the mass)\\}\label{Lem0}
Let $\{v_k\} _{k\in \mathbb{N}}\subset \mathcal{X}$, suppose that there exists a non-negative constant $c$ independent of $k$ such that $E_{\lambda }(v_k)<c$ for all $k\in \mathbb{N}$, then
\begin{equation}\label{eq:compens_mass}
\lim_{l\rightarrow +\infty}\left( \sup_k \left( \int_{l}^{+\infty }\int_{-1}^1 v_k (y,-A)\, dydA\right)\right)=0.
\end{equation}
\end{lem}
\begin{Proofc}{\bf Proof of Lemma \ref{Lem0}\\}
\textbf{Step 1: An estimate for the minimum of $\epsilon _A$}\par
Let $A> 0$ and define
\begin{equation}\label{eq:Compens_1}
m _A =m ^k _A=\int_{-1}^1 v_k (y,-A)\,dy
\end{equation}
Let $\epsilon _A$ be given by \eqref{ene-1D} and $I_{m _A} ^A$ be the minimum of $\epsilon _A$ corresponding to the constraint $m _A$. Then for $\lambda _A$ the root of the trinomial in \eqref{eq:1D_multiplier}, it is
\[
m _A=\frac{(\lambda _A-A)^2(2\lambda _A+A)}{3\lambda _A ^2}.
\]
Using \eqref{form-w} we calculate
\begin{equation}\label{eq:Compens_1.5}
I_{m _A} ^A=\frac{A^2}{a}\left(\frac{2}{3}\frac{(1-a)^3}{a}+(1-a)^2\right) =\frac{(\lambda _A-A)^2(2\lambda _A+A)}{3\lambda _A}=m _A \lambda _A .
\end{equation}
By Step 1 of the proof of Proposition \ref{chara-1D-sol} we have $\lambda _A\geq A+m _A$ hence equation \eqref{eq:Compens_1.5} becomes
\begin{equation}\label{eq:Comens_1.75}
I_{m _A} ^A\geq m _A(A+m _A)\geq Am _A.
\end{equation}
 \textbf{Step 2: The tails of $v_k$ converge uniformly to $0$}\par
We argue by contradiction, suppose that 
\[
\sup_k \left( \int_{l}^{+\infty } m ^k _A \, dA\right) \nrightarrow 0\text{ as } l\rightarrow +\infty
\]
then, there are $\varepsilon >0$ and a sequence $l_j\rightarrow +\infty$ as $j\rightarrow +\infty$ such that 
\begin{equation}\label{eq:Compens_2}
\sup_k \left(\int_{l_j}^{+\infty } m ^k _A \, dA\right) \geq \varepsilon 
\end{equation}
By Fubini's Lemma we have for $l_j>\lambda $
\begin{align}\label{eq:Compens_3}
E_{\lambda }(v_k)&=\int_{0}^{+\infty }\int_{-1}^1\frac{|\nabla v_k|^2}{2}+A|\nabla v_k|-\lambda v_k\,\mathrm{d}y\mathrm{d}A\\
&\geq \int_{0}^{+\infty }\int_{-1}^1 \frac{|\partial _yv_k|^2}{2}+A|\partial _yv_k|-\lambda v_k\, \mathrm{d}y\mathrm{d}A\nonumber\\
&\geq \int_{l_j}^{+\infty } I_{m ^k _A }^A-\lambda m_A ^k\, dA\geq \int_{l_j}^{+\infty } m ^k _A(l_j-\lambda )\, dA, \nonumber
\end{align}
where in the last inequality we used \eqref{eq:Compens_1}, \eqref{min} and \eqref{eq:Comens_1.75}.
Taking the supremum over $k\in \mathbb{N}$ we get, using \eqref{eq:Compens_2}
\[
c\geq \sup_k E_{\lambda }(v_k)\geq \sup_k \left(\int_{l_j}^{+\infty } m ^k _A( l_j-\lambda )\, dA\right) \geq (l_j-\lambda ) \varepsilon \rightarrow +\infty \quad \text{ as } l_j\rightarrow +\infty ,
\]
a contradiction.\\
\end{Proofc}
We have the following Lemma.
\begin{lem}{\bf (Approximation by smooth functions)\\}\label{Lem1}
Let $v\in \mathcal{X}(\hat{\Omega})$. Then, there is a sequence $v_A\in W^{1,2} _0 (\hat{\Omega})$ such that
\begin{equation}\label{eq:Aprox_smooth_1}
v_A\rightarrow v \text{ in } W^{1,2}(\hat{ \Omega})\cap L^1(\hat{ \Omega}) ,
\end{equation}
and 
\begin{equation}\label{eq:Aprox_smooth_2}
\lim_{A\rightarrow +\infty }\int_{\hat{\Omega}} |z||\nabla v_A-\nabla v|= 0.
\end{equation}
\end{lem}
\begin{Proofc}{\bf Proof of Lemma \ref{Lem1}\\}
First we note that $v\in L^1(\hat{\Omega})$ by Remark \ref{Rem:W^1,1}. Let $A>1$ we define the cut off functions $\eta _A\in W^{1,\infty } _0(\mathbb{R})$ by
\[
 \eta _A (z):=\begin{cases} 
1 & \text{ if } |z|\leq A, \\
1-\frac{1}{A}(|z|-A) & \text{ if } A\leq |z|\leq 2A,\\
0 & \text{ if } 2A\leq |z| .
\end{cases}
\]
Then 
\begin{equation}\label{eq:Lem_approx}
|\eta ' _A(z)|\leq \frac{2}{|z|} \text{ a.e. }
\end{equation}
 The functions $v_A(y,z):=\eta _A (z)v(y,z)$ belong to $W^{1,2} (\hat{\Omega})$, they have compact support in $\hat{\Omega} _A=(-1,1)\times (-2A,2A)$ and zero trace on $\partial \hat{\Omega} _A$. Since the boundary of each $\hat{\Omega} _A$ is Lipschitz and bounded we have by \cite[Theorem 15.29]{MR2527916} that $v_A\in W^{1,2} _0 (\hat{\Omega}_A)$. It is not difficult to see that $v_A\rightarrow v$ in $W^{1,2}(\hat{\Omega})$, we will show that $\displaystyle\lim_{A\rightarrow +\infty}|z||\nabla v_A-\nabla v|=0$.\\
 We have 
 \begin{align*}
 \int_{\hat{\Omega}}|z||\nabla v_A-\nabla v|&\leq \int_{\hat{\Omega}}|z|(|\eta '_Av|+|\eta _A-1||\nabla v|)\\
 &\leq \int_{\hat{\Omega} \cap \{A<|z|<2A\}} 2|v|+\int_{\hat{\Omega}\cap \{A<|z|\}}|z||\nabla v|
 \end{align*}
then, using \eqref{eq:Lem_approx} and the fact that $|z||\nabla v|,|v|\in L^1(\hat{\Omega })$ the right hand side of the above estimate converges to zero as $A\rightarrow +\infty$.\\
The convergence in $L^1(\hat{\Omega})$ in \eqref{eq:Aprox_smooth_1} follows by Remark \ref{Rem:W^1,1}.
\end{Proofc}
For two sets $U,U'\subset \mathbb{R}^2$, by $U\subset\subset U'$ we mean that $U$ is relatively compact in $U'$, i.e. $\overline{U}\subset U'$ and $\overline{U}$ is compact. Also for a function $u(y,z)$ we define the positive part $u^+(y,z)=\max \{u(y,z),0\}$.
\begin{Proofc}{\bf Proof of Theorem \ref{theo:main_theo_1} \eqref{itm:E_L_eq_1}\\ }
\textbf{Step 1. Boundedness of $E_\lambda$ from below}\par
We focus in the cases $\lambda >0$ since for $\lambda =0$ the minimizer of $E_\lambda$ is trivially the zero function. We fix $\lambda >0$, let $u\in \mathcal{X}$, using Poincar\'e's inequality in $\Omega$ (Remark \ref{Rem:W^1,1}) we get
\begin{align*}
E_\lambda (u)& =\int_\Omega \frac{|\nabla u|^2}{2}+|z||\nabla u|-\lambda \int_\Omega u\\
&\geq \int_\Omega \frac{|\nabla u|^2}{2}+(|z|-2\lambda )|\nabla u|.
\end{align*}
We split the last integral in the domains $\{|z|\geq 2\lambda \}\cap \Omega$ and $\{|z|\leq 2\lambda \}\cap \Omega$ and get
\begin{align*}
E_\lambda (u)&\geq \int_{\{|z|\leq 2\lambda \}\cap\Omega} \frac{|\nabla u|^2}{2}-2\lambda |\nabla u|\\
&\geq \int_{\{|z|\leq 2\lambda \}\cap\Omega}-2 \lambda ^2>-\infty .
\end{align*}
\textbf{Step 2. Minimizing sequence}\par
Let $u_k\in \mathcal{X}$ with $\displaystyle\lim_{k\rightarrow +\infty }E_\lambda (u_k)= \inf_{v\in\mathcal{X}}E_\lambda (v)$. We will denote by $c$ a generic positive constant which does not depend on the parameter $k$. There is a positive constant $c$ such that $\displaystyle\sup_{k\in \mathbb{N}} E_\lambda (u_k)\leq c$, then as in Step 1 we use Poincare's inequality to get
\begin{align*}
c &\geq \int_{\{|z|\leq 2\lambda \}\cap \Omega} \frac{|\nabla u_k|^2}{2}+(|z|-2\lambda )|\nabla u_k|+\int_{\{|z|\geq 2\lambda \}\cap \Omega}\frac{|\nabla u_k|^2}{2}\\
&\geq \int_{\{|z|\leq 2\lambda \}\cap\Omega }\frac{|\nabla u_k|^2}{2}-\frac{|\nabla u_k|^2}{4}-(|z|-2\lambda )^2+\int_{\{|z|\geq 2\lambda \}\cap\Omega }\frac{|\nabla u_k|^2}{2},
\end{align*}
where in the second inequality we used Young's inequality$\left( |a||b|\leq \frac{b^2}{4} +a^2\right)$. Is is easy now to see that
\begin{equation}\label{eq:dirchlet_unif_bd}
\int_\Omega |\nabla u_k|^2\leq c.
\end{equation}
Then by Poincare's inequality and compactness there is $u\in W^{1,2} _{0L}(\Omega )$ such that $u_k\rightharpoonup u$ as $k\rightarrow +\infty$.

 Using similar arguments we get $c\geq \int_\Omega (|z|-2\lambda )|\nabla u_k|$, or if we split the integral in the domains
 $\{|z|\geq 4\lambda \}\cap \Omega=\{|z|-2\lambda \geq |z|/2\}\cap \Omega$ and $\{|z|\leq 4\lambda \}\cap \Omega$ we get
\begin{equation}\label{eq:unif_bd_totalvar}
\frac{1}{2}\int_{\{|z|-2\lambda \geq |z|/2\}\cap \Omega}|z||\nabla u_k|\leq \int_{\{|z|\geq 4\lambda \}\cap \Omega}(|z|-2\lambda )|\nabla u_k|\leq c-\int_{\{|z|\leq 2\lambda \}\cap \Omega}(|z|-2\lambda )|\nabla u_k| ,
\end{equation}
since $\int_{\{2\lambda \leq |z|\leq 4\lambda \}\cap \Omega}(|z|-2\lambda )|\nabla u_k| \geq 0 $. We can now bound the right hand side of \eqref{eq:unif_bd_totalvar} using H{\"o}lders inequality and \eqref{eq:dirchlet_unif_bd} and get eventually that $\int_{\{|z|\geq  4\lambda \}\cap \Omega}|z||\nabla u_k|\leq c$.
 Using H{\"o}lders inequality and \eqref{eq:dirchlet_unif_bd}, one can also bound the quantity $\int_{\{|z|\leq 4\lambda \}\cap \Omega}|z||\nabla u_k|$ uniformly in $k$, we can therefore conclude that
\begin{equation}\label{eq:final_unif_bd_tot_var}
\int_\Omega |z||\nabla u_k|\leq c,
\end{equation}
where again $c$ is a positive constant independent of $k$.\\
\textbf{Step 3. Lower semicontinuity}\par
We will show that 
\begin{equation}\label{eq:lsc_1}
\int_\Omega \frac{|\nabla u|^2}{2}+|z||\nabla u| \leq \liminf_{k\rightarrow +\infty} \int_\Omega  \frac{|\nabla u_k|^2}{2}+|z||\nabla u_k| ,
\end{equation}
and
\begin{equation}\label{eq:lsc_2}
-\lambda \int_\Omega u\leq  \liminf_{k\rightarrow +\infty }\left( -\lambda \int_\Omega u_k\right)  .
\end{equation}
Equations \eqref{eq:dirchlet_unif_bd}, \eqref{eq:final_unif_bd_tot_var} and \eqref{eq:lsc_1} imply that $u\in \mathcal{X}$ and then $u\in L^1(\Omega )$ by Remark \ref{Rem:W^1,1}. Whereas, equations \eqref{eq:lsc_1} and \eqref{eq:lsc_2} together imply that $E_\lambda (u)\leq \displaystyle\liminf_{k\rightarrow +\infty} E_\lambda (u_k)$, which shows that $u$ is a minimizer of $E_\lambda $ in $\mathcal{X}$. Since the integrand in \eqref{eq:lsc_1} is non-negative convex in the gradient variable and measurable in the $z$ variable, the inequality \eqref{eq:lsc_1} follows from \cite[Chapter I, Theorem 2.5]{MR717034}.

For $l>0$ fixed we have
\begin{align}\label{eq:linear_lsc}
\int_{\Omega } u_k&= \int_{l}^{+\infty }\int_{-1}^1  u_k(y,-A)\,dydA +\int_{0}^{l} \int_{-1}^1 u_k(y,-A)\,dydA\\
&\leq  \sup_k \left(  \int_{l}^{+\infty }\int_{-1}^1  u_k(y,-A)\,dydA \right) +\int_0^{l}\int_{-1}^1  u_k(y,-A)\,dydA.
\end{align}
Since $E_\lambda (u_k)$ is uniformly bounded we
can apply Lemma \ref{Lem0} and get that \eqref{eq:compens_mass} holds for the sequence $u_k$.
Using \eqref{eq:compens_mass} and the fact that $u\in L^1(\Omega )$, we can take the $\limsup$ in \eqref{eq:linear_lsc}, as $k\rightarrow +\infty$ and then $l\rightarrow +\infty$ and get  $\displaystyle\limsup_{k\rightarrow +\infty} \int_\Omega u_k\leq \int_\Omega u$ or else \eqref{eq:lsc_2},
which completes the proof of the lower semi-continuity of $E_\lambda$ and hence the existence of a minimizer $u\in\mathcal{X}$.\\
\textbf{Step 4. Uniqueness}\par
Let $u,\tilde{u}\in \mathcal{X}$ be two minimizers, then using similar arguments as in \cite[Section 3.5.4, p.36]{MR0521262} one can show
that
\begin{equation}\label{eq:uniq_ineq_1}
\int_{\Omega }  \nabla u\cdot (\nabla \tilde{u}-\nabla u)+\int_{\Omega }|z||\nabla \tilde{u}|-\int_{\Omega }|z||\nabla u|\geq \lambda\int_\Omega \tilde{u}-u ,
\end{equation}
\begin{equation}\label{eq:uniq_ineq_2}
\int_{\Omega } \nabla \tilde{u}\cdot (\nabla u-\nabla \tilde{u})+\int_{\Omega }|z||\nabla u|-\int_{\Omega }|z||\nabla \tilde{u}|\geq \lambda \int_\Omega u-\tilde{u} .
\end{equation}
If we add equations \eqref{eq:uniq_ineq_1} and \eqref{eq:uniq_ineq_2} we get
\[
\int_{\Omega }|\nabla u-\nabla \tilde{u}|^2\leq 0,
\]
hence $u=\tilde{u}$ in $\Omega $ since they also have the same lateral boundary conditions.\\
\textbf{Step 5. Non-negative minimizer}\par
We have by \cite[Corollary 2.1.8, page 47]{Ziem} that $\nabla u ^+=(\nabla u)\cdot \chi _{\lbrace u>0\rbrace }$, where by $\chi _{\lbrace u>0\rbrace }$ we denote the characteristic function of the set $\{(y,z):u(y,z)>0\}$. Since also $-\lambda \int_\Omega u^+\leq -\lambda \int_\Omega u$ we have $E_\lambda (u^+)\leq E_\lambda (u)$, hence $u=u^+$ by the uniqueness of minimizers.\\
\textbf{Step 6. $\lambda \in [0,1]$}\par
Our goal is to show that
\begin{equation}\label{eq:Equiv_4}
E_\lambda (u )\geq 0, \quad \text{ for all } u\in \mathcal{X},
\end{equation}
 then because $0\in \mathcal{X}$ and $E_\lambda (0)=0$ we get that the unique minimizer of $E_\lambda$ is the zero function. In view of Lemma \ref{Lem1}, it is enough to prove \eqref{eq:Equiv_4} for functions $u$ with $\hat{u} \in W^{1,2} _0 (\hat{\Omega})$. Let $u$ be such a function, then as in Step 5 we have
\begin{equation}\label{eq:Equiv_5}
E_\lambda (u^+)\leq E_\lambda (u).
\end{equation}
 Suppose that the compact support of $\hat{u}^+$ is contained in $[-1,1]\times (-A,A)$ where $A$ is large enough, then we have
\begin{align*}
\int_{\Omega}|z||\nabla u^+|&\geq \int_{-1}^1\int_{0}^A|z|\left|  \frac{\partial u^+}{\partial z}\right|   \, dzdy\\
&\geq \left| \int_{-1}^1\int_{0}^A z\frac{\partial u^+}{\partial z}\, dzdy\right|  \\
& =\int_{\Omega}u^+
\end{align*}
in the last equality we used integration by parts. This estimate together with \eqref{eq:Equiv_5} and the fact that $\lambda \leq 1$ gives
\[
E_\lambda (u)\geq E_\lambda (u^+)\geq (1-\lambda )\int_{\Omega}u^+\geq 0.
\]
\textbf{Step 7. $\lambda \in ( 1,+\infty )$}\par
Our goal is to prove that there is $u\in \mathcal{X}$ with $E_\lambda (u)<0$. Let $\varphi \in C^\infty (-1,1)$, $\varphi \geq 0$ with $\varphi (-1)=0=\varphi (1)$ and $\int_{-1}^1\varphi =1$ (for example $\varphi (y)=\frac{3}{4}(1-y^2)$). We define
\[
u(y,z):=k^{-3}e^{kz} \varphi (y)
\]
where $k>0$ is large enough, to be chosen later. It is $u\in \mathcal{X}$ and 
\begin{align}\label{eq:Equiv_6}
\int_\Omega \frac{|\nabla u|^2}{2}&=\frac{1}{2}\int_{-1}^1 \left[ k^{-2}(\varphi '(y))^2 +\varphi ^2 (y)  \right]dy \int_{-\infty }^0 (k^{-2}e^{kz})^2\, dz \\ &=\frac{A_k}{4}k^{-5}\nonumber
\end{align}
where we set $A_k=\int_{-1}^1 \left[ k^{-2}(\varphi '(y))^2 +\varphi ^2 (y)  \right]\, dy$. Also
\begin{align}\label{eq:Equiv_7}
\int_\Omega |z||\nabla u|&= \int_{-1}^1 \sqrt{k^{-2}(\varphi '(y))^2 +\varphi ^2 (y)}\, dy \int_{-\infty }^0 |z|k^{-2}e^{kz}\, dz \\
&=B_k\int_{-\infty }^0 |z|k^{-2}e^{kz}\, dz\nonumber
\end{align}
where $B_k=\int_{-1}^1 \sqrt{k^{-2}(\varphi '(y))^2 +\varphi ^2 (y)}\, dy$. If we integrate by parts the second product component of the right hand side of \eqref{eq:Equiv_7} we get
\[
\int_{-\infty }^0 |z|k^{-2}e^{kz}\, dz=\int_{-\infty }^0 k^{-3}e^{kz}=\int_{-\infty }^0 k^{-3}e^{kz}\int_{-1}^1\varphi =\int_\Omega u,
\]
then \eqref{eq:Equiv_7} becomes
\begin{equation}\label{eq:Equiv_8}
\int_\Omega |z||\nabla u|=B_k\int_\Omega u.
\end{equation}
We also have $\int_\Omega u=k^{-4}$, then we can write $E_\lambda (u)$ using \eqref{eq:Equiv_6} and \eqref{eq:Equiv_8} as
\begin{equation}\label{eq:Equiv_9}
E_\lambda (u)= \frac{A_k}{4}k^{-5}+(B_k-\lambda )k^{-4}.
\end{equation}
Next we note that $B_k\geq 1$, is decreasing in $k$ (and so is $A_k$) and $B_k\rightarrow 1$ as $k\rightarrow +\infty$. Since $\lambda >1$ we can find $k_0$ large enough such that $B_{k_0}<\lambda $, then \eqref{eq:Equiv_9} becomes
\[
E_\lambda (u)\leq  \frac{A_{k_0}}{4}k^{-5}+ (B_{k_0}-\lambda )k^{-4},
\]
for all $k\geq k_0$. We can now conclude if we choose $k\geq k_0$ large enough, since the function $k^{-5}$ decreases faster than $k^{-4}$, for example $k>\max\{ k_0, \frac{A_{k_0}}{4(\lambda -B_{k_0})}\}$.
\end{Proofc}
\subsection{The $\varepsilon$-approximation}\label{Subsec:Regul}
Let $\lambda >0$, $u_\lambda$ be the minimizer of $E_\lambda $ given by Theorem \ref{theo:main_theo_1} \eqref{itm:E_L_eq_1}. For $A>0$ we define $\hat{\Omega} _A=\{(y,z)\in \hat{\Omega}:\, |z|\leq A\}$, $\Omega _A=\Omega \cap\hat{\Omega} _A$ and
\[
\mathcal{H}_A=\lbrace v\in W^{1,2} (\Omega _A), \, v=u_\lambda , \text{ on } \partial \Omega _A \setminus \lbrace z=0\rbrace \rbrace ,
\]
We are interested in approximate minimizers of \eqref{eq:energy_problem}, for this we study the minimizers in $\mathcal{H}_A$ of the approximate functional
\begin{equation}\label{eq:approx_energy}
E^A _{\varepsilon ,\lambda } (u)= \int_{\Omega _A} \frac{|\nabla u |^2}{2}+|z|\sqrt{\varepsilon ^2+|\nabla u |^2}-\lambda u,
\end{equation}
where $\varepsilon >0$.\\
Since we have mixed boundary conditions, an easy way to describe the space of test functions for the weak formulation of the first variation of \eqref{eq:approx_energy} is to use reflection in the domain $\hat{\Omega} _A$. We will simply write $\hat{\phi}\in W^{1,2} _0({\hat{\Omega}_A})$ for the test functions.
We have the following Proposition.
\begin{pro}{\bf ($W^{2,2}_{\mathrm{loc}}$ regularity of approximate problem)\\}\label{Pro:Reg_approx}
Let $A, \varepsilon ,\lambda >0$, then there exists a unique minimizer $u_{\varepsilon ,A}\in \mathcal{H} _A $ of $E^A _{\varepsilon ,\lambda}$. Moreover, $\hat{u}_{\varepsilon ,A} \in W^{2,2}_{\mathrm{loc}}(\hat{\Omega}_A)$ and the following equation holds
\begin{equation}\label{eq:E_L_appr}
\int_{\Omega _A } \nabla u_{\varepsilon ,A}\cdot \nabla \varphi +|z|\frac{ \nabla u_{\varepsilon ,A}\cdot \nabla \varphi}{\sqrt{\varepsilon ^2+|\nabla u_{\varepsilon ,A}|^2}} =\lambda \int_{\Omega _A} \varphi ,\quad \text{ for all } \hat{\varphi}\in W^{1,2} _0(\Omega _A) ,
\end{equation} 
and $\partial _z u_{\varepsilon ,A} (y,0)=0$ for $y\in (-1,1)$.
\end{pro}
The existence of a minimizer is a consequence of the direct method in the bounded domain $\Omega _A$, while the regularity results are standard. We give a sketch of the Proof of Proposition \ref{Pro:Reg_approx} in Appendix \ref{appendix}.
\begin{Proofc}{\bf Proof of Theorem \ref{theo:main_theo_1} \eqref{itm:E_L_eq_2}-\eqref{itm:E_L_eq_4}\\}
\textbf{Step 1. Solutions of E-L equation are minimizers of \eqref{eq:energy_problem}}\par
First we will show that for any pair $(u,q)\in \mathcal{X}\times\Lambda$ that satisfies equation \eqref{eq:sub_differential}, $u$ is a minimizer of $E_\lambda $. Let $v\in \mathcal{X}$, using \eqref{eq:sub_differential} and the fact that $|q|\leq 1$ it is easy to check that $q\in \partial |\cdot |(\nabla u)$ in $\Omega$. By the definition of the subdifferential we have
\begin{equation}\label{eq:E_L__eq_1}
E_\lambda (v)-E_\lambda (u)\geq  \int_\Omega \nabla u\cdot \nabla (v-u) +|z|q\cdot \nabla (v-u)-\lambda \int_\Omega (v-u)=0,
\end{equation}
where we used \eqref{eq:sub_differential} with test function $\varphi =v-u\in \mathcal{X}$.\\
\textbf{Step 2. Approximating solutions}\par
As usual we will focus in the case $\lambda >0$. Let $u=u_\lambda$ be the minimizer of $E_\lambda$ given by Theorem \ref{theo:main_theo_1} \eqref{itm:E_L_eq_1}. For $\varepsilon >0$ let $u_{\varepsilon ,A}$ be the minimizer of $E^A _{\varepsilon ,\lambda}$ given by Proposition \ref{Pro:Reg_approx}, then for all $A>0$ we will show that $u_{\varepsilon ,A} \rightarrow u$ strongly as $\varepsilon \rightarrow 0$, in $W^{1,2}(\Omega _A)$ up to a subsequence. Extending $u_{\varepsilon ,A}$ by $u_\lambda$ outside $\Omega _A$, we can write the following variational inequalities as in the Step 1 of the proof of Theorem \ref{theo:main_theo_1} \eqref{itm:E_L_eq_1}
\begin{equation}\label{eq:appr_var_ineq_1}
\int_{\Omega _A}\nabla u\cdot (\nabla u_{\varepsilon ,A} -\nabla u)+\int_{\Omega _A}|z||\nabla u_{\varepsilon ,A} |-\int_{\Omega _A}|z||\nabla u|\geq \lambda\int_{\Omega _A} u_{\varepsilon ,A} -u
\end{equation}
and
\begin{equation}\label{eq:appr_var_ineq_2}
\int_{\Omega _A}\nabla u_{\varepsilon ,A} \cdot (\nabla u-\nabla u_{\varepsilon ,A} )+\int_{\Omega _A}|z|\sqrt{\varepsilon ^2+|\nabla u|^2}-\int_{\Omega _A}|z|\sqrt{\varepsilon ^2+|\nabla u_{\varepsilon ,A}|^2}\geq \lambda \int_{\Omega _A} u -u_{\varepsilon ,A}  .
\end{equation}
Adding inequalities \eqref{eq:appr_var_ineq_1} and \eqref{eq:appr_var_ineq_2}, we get
\begin{align*}
\int_{\Omega _A}|\nabla u_{\varepsilon ,A} -\nabla u|^2&\leq \int_{\Omega _A}|z|(|\nabla u_{\varepsilon ,A} |-\sqrt{\varepsilon ^2+|\nabla u_{\varepsilon ,A} |^2})  +|z|(\sqrt{\varepsilon ^2+|\nabla u |^2}-|\nabla u|)\\
&\leq \int_{\Omega _A}|z|(\sqrt{\varepsilon ^2+|\nabla u|^2}-|\nabla u|)\\
&=\int_{\Omega _A}|z|\frac{\varepsilon^2}{\sqrt{\varepsilon ^2+|\nabla u|^2}+|\nabla u|}\leq A|\Omega _A|\varepsilon ,
\end{align*}
Then using also Poincare's inequality we get for all $A>0$ and up to a subsequence
\begin{equation}\label{eq:der_conv}
\nabla u_{\varepsilon ,A} \rightarrow \nabla u,\,u_{\varepsilon ,A} \rightarrow u \text{ a.e.  in }\Omega _A \text{ as }\varepsilon \rightarrow 0.
\end{equation}
\textbf{Step 3. The function $q$}\par
For $q_{\varepsilon ,A} =\frac{\nabla u_{\varepsilon ,A}}{\sqrt{ \varepsilon ^2+|\nabla u_{\varepsilon ,A}|^2}}$, we have $q_{\varepsilon ,A}\cdot \nabla u_{\varepsilon ,A} \leq |\nabla u_{\varepsilon ,A}|$, then using \eqref{eq:der_conv} it is not difficult to see that 
\begin{equation}\label{eq:norm_strong_conv}
\displaystyle q_{\varepsilon ,A}\cdot \nabla u_{\varepsilon ,A} \rightarrow |\nabla u| \text{ a.e. in } \Omega _A \text{ as } \varepsilon \rightarrow 0.
\end{equation}
Since $q_{\varepsilon ,A}\in L_{\mathrm{loc}} ^2(\Omega _A ,\mathbb{R}^2)$ with $|q_{\varepsilon ,A}|\leq 1$, there exists $q_A\in L_{\mathrm{loc}} ^2(\Omega _A ,\mathbb{R}^2)$ with $|q_{A}|\leq 1$ and such that $q_{\varepsilon ,A}$ converges weakly to $q_A$ in $L^2 (U,\mathbb{R}^2)$, as $\varepsilon \rightarrow 0$, for every $U\subset\subset \Omega _A$. Then using also \eqref{eq:der_conv} we have $\displaystyle\lim_{\varepsilon \rightarrow 0} \int_U q_{\varepsilon ,A}\cdot \nabla u_{\varepsilon ,A} =\int_U q_A\cdot \nabla u $ for all $U\subset \subset \Omega _A$ and by \eqref{eq:norm_strong_conv} we get that $q_A\cdot \nabla u=|\nabla u|$ a.e. in $\Omega _A$. Extending $q_A$ by zero outside $\Omega _A$ we may wright $q_A\in L^2 _{\mathrm{loc}}(\Omega ,\mathbb{R}^2)$ and as before we can find $q\in L^2 _{\mathrm{loc}}(\Omega ,\mathbb{R}^2)$, with $|q|\leq 1$ and such that $q_A$ converges weakly to $q$ in $L^2 (U,\mathbb{R}^2)$, as $A \rightarrow +\infty$, for every $U\subset\subset \Omega _A$, and hence $q\cdot \nabla u=|\nabla u|$ a.e.
\\
\textbf{Step 4. Passing to the limit $\varepsilon \rightarrow 0, A\rightarrow +\infty$}\par
Let $\varphi$ with $\hat{\varphi}\in W^{1,2} _0(\hat{\Omega})$, then equation \eqref{eq:E_L_appr} with $A$ large enough holds for this test function and since $q_{\varepsilon ,A}$ is bounded we can pass to the limit as $\varepsilon\rightarrow 0$ and get
\[
\int_{\Omega} \nabla u\cdot \nabla\varphi +|z|q_A\cdot \nabla\varphi =\lambda \int_{\Omega} \varphi .
\]
We can now pass to the limit as $A\rightarrow +\infty$ and using also Lemma \ref{Lem1} we get \eqref{eq:sub_differential}.\\
\textbf{Step 5. Uniqueness}\par
Let $(u,q )$, $(\bar{u},\bar{q})$ be two solutions of \eqref{eq:sub_differential} then by Step 1 we have $u=\bar{u}$, since minimizers of \eqref{eq:energy_problem} in $\mathcal{X}$ are unique by Theorem \ref{theo:main_theo_1} \eqref{itm:E_L_eq_1}. Then in the set $\{\nabla u\neq 0\}$ the vectors $q,\bar{q}$ are parallel to $\nabla u$ and so is $q-\bar{q}$, but since $(q-\bar{q})\cdot \nabla u=0$ by \eqref{eq:sub_differential} we have $q=\bar{q}$ a.e. in $\{\nabla u\neq 0\}$.\\
\textbf{Step 6. Neumann condition}\par
We denote by $\partial _{x_i} $, $i=1,2$ respectively the derivatives $\partial _y,\, \partial _z$. Let $i,j\in \{1,2\}$, $\hat{U}\subset\subset \hat{\Omega}_A$, by Proposition \ref{Pro:Reg_approx} we have that $\hat{u} _{\varepsilon ,A} \in W^{2,2} _{\mathrm{loc}}(\hat{U} )$, by Lemma \ref{lem:unif_bd_2der} the second derivatives of $\hat{u} _{\varepsilon ,A}$ are uniformly bounded in $L^2(\hat{U})$, hence for $\varphi \in W^{1,2} _0(\hat{U})$ we have (up to a subsequence)
\[
\int_{\hat{U}} \partial _{x_i}\hat{u} \partial _{x_j}\varphi =\lim_{\varepsilon \rightarrow 0}\int_{\hat{U}}\partial _{x_i}\hat{u} _{\varepsilon ,A}\partial _{x_j}\varphi=-\lim_{\varepsilon \rightarrow 0}\int_{\hat{U}}\partial _{x_j}\partial _{x_i}\hat{u} _{\varepsilon ,A} \varphi=-\int_{\hat{U}} g\varphi ,
\]
for some function $g\in L^2(\hat{U})$. We have proved that $\hat{u} \in W^{2,2} _{\mathrm{loc}}(\hat{\Omega} )$, then applying a Sobolev embedding Theorem (\cite[Section 5.6.3]{Eva}) we get that $\hat{u} \in C^{0,\alpha } _{\mathrm{loc}}(\hat{\Omega} )$ for all $\alpha \in (0,1)$. 
As in the proof of Proposition \ref{Pro:Reg_approx} we can now define the trace of the derivative of $u$ on $\{z=0\}$ and $\partial _zu(y,0)=0$ for $y\in (-1,1)$.
\end{Proofc}
\section{Properties of the solution}\label{SecProperties}
\subsection{Comparison Principle}\label{SecComparison}
In view of Theorem \ref{theo:main_theo_1} \eqref{itm:E_L_eq_1} we will assume that $\lambda >1$ for the rest of the paper.
\begin{defi}{\bf Sub/supersolution\\}
Let $u\in \mathcal{X}$ be non-negative and $q\in \Lambda$, $\Lambda$ as in \eqref{eq:lambda_def}, we call the pair $(u,q)$ a subsolution (resp. a supersolution) of the equation \eqref{eq:sub_differential} if
\begin{equation}\label{eq:defi_subsol}
\begin{cases} 
\int_\Omega \nabla u\cdot \nabla \varphi +|z|q\cdot \nabla \varphi \leq \lambda \int_\Omega \varphi \,\,(\text{resp.} \geq \lambda \int_\Omega \varphi )& \text{for all } \varphi \in \mathcal{X}, \varphi \geq 0,\\
q\cdot \nabla u=|\nabla u|  &  \text{a.e. in }\Omega .
\end{cases} 
\end{equation}
\end{defi}
\begin{pro}\label{comp-prin}{\bf Comparison principle}\\
Let $u,v\in \mathcal{X}$, $q_u,q_v\in \Lambda$ with $(u,q_u),(v,q_v)$ a subsolution and a supersolution respectively of \eqref{eq:sub_differential}, with $0=u\leq v$ on $\{-1,1\}\times (-\infty ,0)$ in the sense of traces, then
\[
u\leq v,\quad \text{ in } \Omega .
\]
\end{pro} 
\begin{Proofc}{\bf Proof of Proposition \ref{comp-prin}}\\
Let $\varphi =(u-v)_+$, then $\varphi\in \mathcal{X}$. If we write the inequalities \eqref{eq:defi_subsol} for $u,v$ with this test function and subtract the one from the other we get
 \[
 \int_\Omega \nabla (u-v)\cdot \nabla (u-v)_++|z|(q_u-q_v)\cdot \nabla (u-v)_+\leq 0,
 \]
 or if we use \cite[Corollary 2.1.8, page 47]{Ziem} we can write it as
 \begin{equation}\label{eq:compar_ineq1}
 \int_\Omega |\nabla (u-v)|^2\chi _{\{u-v\geq 0\}}\leq -\int_\Omega |z|\left[ (q_u-q_v)\cdot \nabla (u-v)\right] \chi _{\{u-v\geq 0\}} .
 \end{equation}
 Next we calculate, using  the properties of $q_u,q_v$ in Definition \ref{eq:defi_subsol}
 \begin{align*}
 (q_u-q_v)\cdot (\nabla u-\nabla v)& =|\nabla u|-q_u\cdot \nabla v-q_v\cdot \nabla u+|\nabla v|\\
 &\geq |\nabla u|-|\nabla u|+|\nabla v|-|\nabla v|=0,\, \text{ a.e.}
 \end{align*}
 then \eqref{eq:compar_ineq1} implies 
 \[
 \nabla (u-v)=0,\quad \text{ a.e. in }\{u-v\geq 0\} 
 \]
 or $\nabla (u-v)_+=0$ almost everywhere. Using the boundary conditions we can conclude that $(u-v)_+=0$ and hence $u\leq v$ a.e. in $\Omega$.
\end{Proofc}
\begin{rem}{\bf (Monotonicity in $\lambda$)\\}
For $u_\lambda$ the minimizer of $E_\lambda $ in $\mathcal{X}$ and $m(\lambda )=\int_\Omega u_\lambda$ the volume rate, using the comparison principle from Proposition \ref{comp-prin} it is not difficult to see that $m(\lambda )$ is increasing in $\lambda$. Unfortunately, the physical volume rate is given, using the rescaling \eqref{eq:change_mass}, by $m_0=(l^2\mu _s g_0\cos\theta )m$, which does not allow us to directly study the monotonicity with respect the inclination angle $\theta$ ($\cos\theta$ is decreasing for $\theta \in [0,\pi /2)$ and $\lambda (\theta )$ is increasing by \eqref{eq:abs_rel_lagr_mult}).
\end{rem}

\subsection{Some explicit profiles}
\begin{figure}[htbp]
     \centering
     \subfloat[][$\phi _{K(\lambda )}$ for $\lambda =1.2,\, 1.4,\, 1.6,\, 1.8.$]{\includegraphics[scale=0.8]{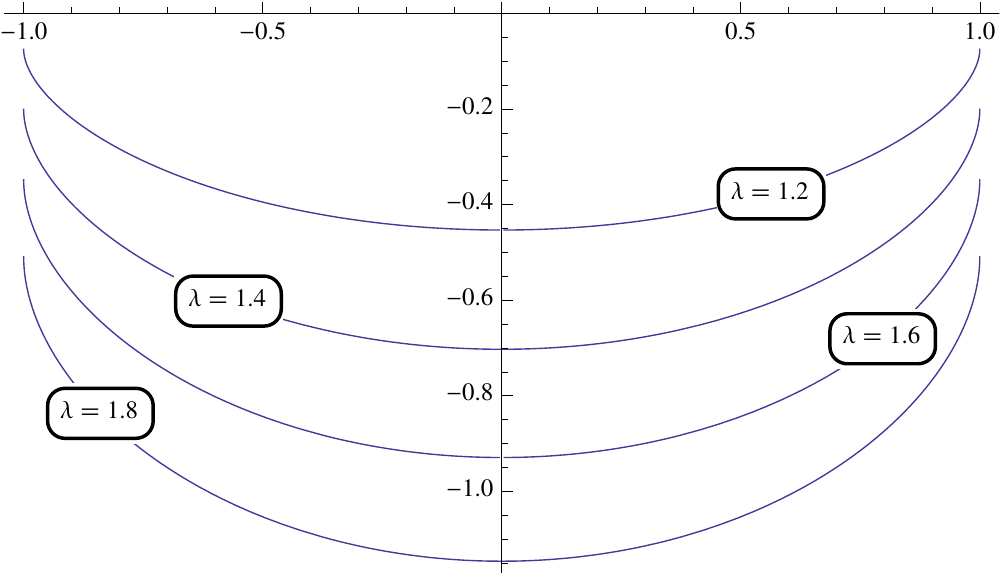}\label{fig:lambdas}}
     \subfloat[][$\phi _{K(\lambda )}(1)$]{\includegraphics[scale=0.8]{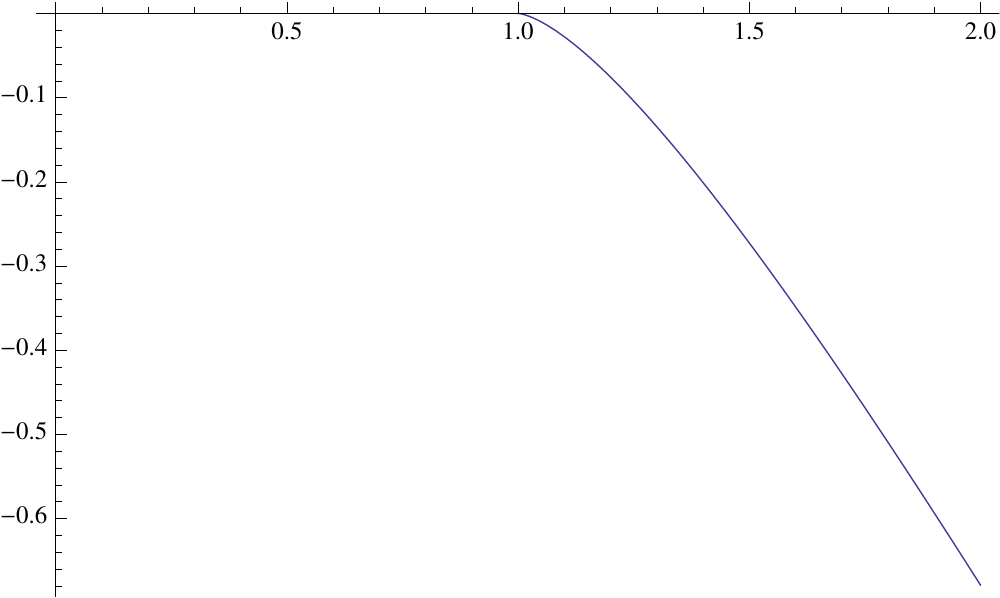}\label{fig:decreasing}}
     \caption{ }
     \label{fig:avalance_2}
\end{figure}

As we explained in the introduction, we study the first variation of the functional \eqref{eq:one_dim_cur}, i.e.
  \begin{equation}\label{eq:eq:E-L_explic_energy}
\frac{\phi \phi ''}{(1+|\phi '|^2)^{3/2}}-\frac{1}{\sqrt{1+|\phi '|^2}}+\lambda =0,\quad y\in (-1,1).
\end{equation}
\begin{lem}{\bf(An explicit solution of \ref{eq:eq:E-L_explic_energy})\\}\label{Lem:explicit_minim}
Let $\lambda >1$, $K(\lambda )$ be given by \eqref{eq:K_max} and $\phi _{K(\lambda )}$ defined in \eqref{eq:expl_implicite_solution}. Then the function $\phi _{K(\lambda )}\in C^\infty (-1,1)\cap C([-1,1])$  is non-positive and the following properties hold
\begin{equation}\label{eq:proper_explicit_minim}
\displaystyle\lim_{y\rightarrow -1}\phi _{K(\lambda )}'(y)=-\infty ,\, \displaystyle\lim_{y\rightarrow 1}\phi _{K(\lambda )}'(y)=+\infty .
\end{equation}
Moreover the function $\phi _{K(\lambda )}$ is convex with minimum $\phi _{K(\lambda )}(0)=\frac{K(\lambda )}{\lambda -1}$ and maximum $\phi _{K(\lambda )}(\pm 1)=\frac{K(\lambda )}{\lambda}$ and if $\bar{\lambda}>\lambda$ then $\phi _{K(\bar{\lambda})}(y)<\phi _{K(\lambda )}(y)$, for $y\in [-1,1]$.
\end{lem}
\begin{Proofc}{ \bf Proof of Lemma \ref{Lem:explicit_minim}\\ }
\textbf{Step 1. The inverse function}\par
Let $\lambda >1$ and $Z\in [\frac{1}{\lambda},\frac{1}{\lambda -1}]$, $f_\lambda (Z)$ be given by \eqref{eq:inverce_funct}. Notice that $f_\lambda$ is smooth in $(\frac{1}{\lambda},\frac{1}{\lambda -1})$ and that it has been chosen so that
\begin{equation}\label{eq:invers_deriv}
f' _\lambda (Z)=\frac{(\lambda Z -1)\sqrt{\lambda ^2-1}}{\sqrt{1-((\lambda ^2-1)Z-\lambda )^2 } },
\end{equation} 
from which we get that $f_\lambda$ is strictly increasing in $[\frac{1}{\lambda},\frac{1}{\lambda -1}]$. We set
\begin{equation}\label{eq:A_lambda}
A_\lambda := f_\lambda \left( \frac{1}{\lambda -1}\right)-f_\lambda \left(\frac{1}{\lambda}\right) = \frac{\pi}{2(\lambda ^2-1)^{3/2}}+\frac{1}{\lambda ^2-1}\left(  1+\frac{\mathrm{Arcsin}\left( \frac{1}{\lambda }\right)}{\sqrt{\lambda ^2-1}} \right),
\end{equation}
by the monotonicity of $f$ we can define the positive function $\phi$ implicitly in the intervals $[-A_\lambda ,0]$ and $[0,A_\lambda ]$ as follows
\begin{equation}\label{eq:expl_implicite_solution_1}
f_\lambda (\phi (y))=f_\lambda \left( \frac{1}{\lambda -1}\right) -|y| ,\quad y\in [-A_\lambda ,A_\lambda ],
\end{equation}
then $f_\lambda (\phi (y))=f_\lambda (\phi (-y))$ for $y\in [0,A_\lambda ]$, which means that $\phi$ is an even function thanks to the monotonicity of $f_\lambda$.
Also by \eqref{eq:expl_implicite_solution_1} we have $\phi (0)=1/(\lambda -1)$
and by \eqref{eq:invers_deriv} we can calculate the limit $\displaystyle\lim_{Z\rightarrow 1/(\lambda -1)} f'(Z)$ and get $\displaystyle\lim_{y\rightarrow 0^+}\phi '(y)=0$. Since $\phi$ is even and smooth in the intervals $[-A_\lambda ,0)$ and $(0,A_\lambda ]$ we eventually get $\phi '(0)=0$.
We have concluded that $\phi \in C^1 (-A_\lambda  ,A_\lambda )$.

 Relation \eqref{eq:expl_implicite_solution_1} gives also for $y\in [-A_\lambda ,A_\lambda ]$
\begin{equation}\label{eq:value_range}
1/\lambda =\phi (\pm A_\lambda )\leq \phi (y)\leq \phi (0)= 1/(\lambda -1)
\end{equation}
 and by \eqref{eq:invers_deriv}
\begin{equation}\label{eq:Lem_expl_min_neumann}
\phi '(-A_\lambda )=+\infty ,\, \phi '(A_\lambda )=-\infty .
\end{equation}
\textbf{Step 2. $\phi$ satisfies \eqref{eq:eq:E-L_explic_energy}}\par
Using \eqref{eq:invers_deriv} we can differentiate \eqref{eq:expl_implicite_solution_1} and taking the squares in both sides of the equation, we get for $y\in (-A_\lambda ,A_\lambda ),$
\[
|\phi '|^2 \frac{(\lambda \phi -1)^2(\lambda ^2-1)}{1-((\lambda ^2-1)\phi -\lambda )^2}=1
\]
or after a few simplifications
\[
|\phi '|^2=\frac{1}{(\lambda -\frac{1}{\phi})^2}-1 .
\]
Noting that $\phi \geq 1/\lambda >0$, the above equation can be rewritten as
\begin{equation}\label{eq:Lem_explici_min_diff_eq}
\phi \left( \lambda -\frac{1}{\sqrt{1+|\phi '|^2}}\right) =1.
\end{equation}
Let $K_0<0$, we define
\begin{equation}\label{eq:1_param_function}
\phi _{K_0}(y):=K_0\phi (\frac{y}{K_0}),\quad y\in [ A_\lambda K_0,-A_\lambda K_0],
\end{equation}
by \eqref{eq:Lem_explici_min_diff_eq}, the negative function $\phi _{K_0}$ satisfies
\begin{equation}\label{eq:Lem_expl_min_diff_eq_2}
\phi _{K_0}(y)\left( \lambda -\frac{1}{\sqrt{1+|\phi _{K_0} '(y)|^2}}\right) =K_0,\quad y\in (A_\lambda K_0,-A_\lambda K_0) .
\end{equation}
In particular, if $K(\lambda )$ is given by \eqref{eq:K_max}, differentiating \eqref{eq:Lem_expl_min_diff_eq_2} with respect to $y$ we get
\begin{equation}\label{eq:eq_of_phi_K_max}
\phi _{K(\lambda )} '\left( \frac{\phi _{K(\lambda )}\phi _{K(\lambda )} ''}{(1+|\phi _{K(\lambda )} '|^2)^{3/2}}-\frac{1}{\sqrt{1+|\phi _{K(\lambda )}'|^2}}+\lambda \right)=0,\quad y\in (-1,0)\cup (0,1).
\end{equation}
 Using equation \eqref{eq:eq_of_phi_K_max} we calculate for $y\in (-1,0)\cup (0,1)$
\begin{equation}\label{eq:second_der}
\phi _{K(\lambda )} ''=\frac{(1+|\phi _{K(\lambda )} '|^2)(\lambda \sqrt{1+|\phi _{K(\lambda )} '|^2}-1)}{-\phi _{K(\lambda )}} > 0,
\end{equation}
here we have also used equation \eqref{eq:Lem_expl_min_diff_eq_2} in order to get the sign of the second derivative. Since $\phi _{K(\lambda )}\in C^1(-1,1)$ we get from \eqref{eq:second_der} that in fact $\phi _{K(\lambda )}\in C^2((-1,1))$. Differentiating further \eqref{eq:second_der} and using \eqref{eq:value_range} we get by iteration $\phi _{K(\lambda )}\in C^\infty (-1,1)$. \\
\textbf{Step 3. Extrema}\par
 By \eqref{eq:value_range} and \eqref{eq:1_param_function} we have
\begin{equation}\label{eq19}
\left\{
\begin{aligned}
&\min_{|y|\leq 1}\phi _{K(\lambda )}=\phi  _{K(\lambda )}(0)=\frac{K(\lambda )}{\lambda -1}=\frac{1}{(\lambda -1)(f_\lambda (\frac{1}{\lambda})-f_\lambda (\frac{1}{\lambda -1}))} ,\\
&\max_{|y|\leq 1}\phi _{K(\lambda )}=\phi _{K(\lambda )}(-1)=\phi _{K(\lambda )}(1)=\frac{K(\lambda )}{\lambda }=\frac{1}{\lambda (f_\lambda (\frac{1}{\lambda})-f_\lambda (\frac{1}{\lambda -1}))}.
\end{aligned}
\right.
\end{equation}
It is
\begin{equation}\label{eq:lem_expl_min_3}
\frac{\mathrm{d}}{\mathrm{d}\lambda }\phi _{K(\lambda )}(1)=-\frac{4(2\lambda ^2+1)\mathrm{Arcsin}\left( \frac{1}{\lambda} \right)\sqrt{\lambda ^2-1} +2\pi (2\lambda ^2+1)\sqrt{\lambda ^2-1}+4(\lambda ^2-1)(\lambda ^2+2)}{\lambda ^2\left(  2\sqrt{\lambda ^2-1}+\pi +2\mathrm{Arcsin}\left( \frac{1}{\lambda } \right)   \right) ^2} <0
\end{equation}
and
\begin{equation}\label{eq:lem_expl_min_4}
\frac{\mathrm{d}}{\mathrm{d}\lambda }\phi _{K(\lambda )}(0)=-\frac{4\left( \lambda -\frac{1}{2}\right) (\lambda -1 )\sqrt{\lambda ^2-1}\left( \pi +2 \mathrm{Arcsin}\left( \frac{1}{\lambda} \right)\right) +4(\lambda ^2-1)\left[ (\lambda -1 )^2+\frac{\lambda -1}{\lambda}\right]       }{(\lambda -1)^2\left(  2\sqrt{\lambda ^2-1}+\pi +2\mathrm{Arcsin}\left( \frac{1}{\lambda } \right)   \right) ^2}<0.
\end{equation}
Figure \ref{fig:decreasing} is the graph of the function $\phi _{K(\lambda )}(1)$ in terms of the variable $\lambda .$\\
\textbf{Step 4. Monotonicity of the graphs in $\lambda$}\par
Let $\bar{\lambda}>\lambda$ we will show that $\phi _{K(\bar{\lambda} )}(y)<\phi _{K(\lambda )} (y)$, for $y\in [0,1]$. Since the functions are even and we already have the monotonicity of the boundary points by  Step 3, we will focus in the interval $(0,1)$. If we use equation \eqref{eq:eq_of_phi_K_max}, we get that the function $w(y)=\phi _{K(\bar{\lambda} )}(y)-\phi _{K(\lambda )} (y)$ satisfies the elliptic equation
\[
-a_1(y)w''(y)+a_2(y)w'(y)+a_3(y)w(y)=\lambda -\bar{\lambda} ,
\]
with 
\[
a_1(y)=\frac{-\phi _{K(\bar{\lambda})}(y)}{(1+|\phi '_{K(\bar{\lambda})}(y)|^2)^{3/2}}, \quad a_3(y)=\frac{\phi ''_{K(\lambda )}(y)}{(1+|\phi '_{K(\bar{\lambda})}|^2)^{3/2}},
\]
and
\[
a_2(y)=\int_0^1 G_1(p(t,y))dt+\phi ''_{K(\lambda )}(y)\phi _{K(\lambda )}(y)\int_0^1 G_2(p(t,y))dt,
\]
with $p(t,y)=\phi '_{K(\lambda )}(y)+t(\phi '_{K(\bar{\lambda})}(y)-\phi '_{K(\lambda )}(y))$, $G_1(p)=\frac{-p}{(1+|p|^2)^{3/2}}$ and $G_2(p)=\frac{-3p}{(1+|p|^2)^{5/2}}$. It is $a_i\in C(0,1)$, $i=1,2,3$ with $a_1,a_3>0$ in $(0,1)$ and $w\in C^2((0,1))\cap C([0,1])$ with $w(0),w(1)<0$ by \eqref{eq:lem_expl_min_3}, \eqref{eq:lem_expl_min_4}. We can now conclude that $w<0$ by a maximum principle.
\end{Proofc}
Using the function $\phi _{K(\lambda )}$ constructed in Lemma \ref{Lem:explicit_minim} we can define a diffeomorphism in $\mathcal{C}_\lambda \cap \Omega $, with $\mathcal{C}_\lambda$ as in \eqref{eq:cone}. Let $L\in (0, +\infty )$, we define 
\[
\phi _L(y):=L\phi _{K(\lambda )}\left( \frac{y}{L} \right) ,\quad y\in [-L,L] .
\]
 We have the following Lemma.
\begin{lem}{\bf(A diffeomorphism)\\}\label{Lem_subsol_2}
Let $\phi _{K(\lambda )}$ be as in \eqref{eq:expl_implicite_solution}, then for $(y,z)\in \overline{\mathcal{C}_\lambda \cap \Omega}\setminus \{(0,0)\} $ there is a unique $L=L(y,z)\in (0,+\infty )$ implicitly defined by 
    \begin{equation}\label{eq:implicit_def_K}
    z=L\phi _{K(\lambda )} \left( \frac{y}{L}\right)=\phi _L (y) ,
    \end{equation}
   and $L\in C^\infty (\mathcal{C}_\lambda \cap \Omega )\cap C(\overline{\mathcal{C}_\lambda \cap \Omega}\setminus \{(0,0)\})$.
\end{lem}
\begin{Proofc}{\bf Proof of Lemma \ref{Lem_subsol_2}\\ }
Since the family of curves $\{(y,\phi _L (y))\}_{L\in (0,+\infty )}$ are obtained as a rescaling of the function $\phi _{K(\lambda )}$ we have that the mapping $(y,L)\mapsto (y,z)$ is a surjection; it is also an injection since the family of curves $\{(y,\phi _L (y))\}_{L\in (0,+\infty )}$ do not intersect. On the other hand the same bijective correspondence can be established locally by the implicit function theorem since $\frac{y}{L}\phi _{K(\lambda )} '\left( \frac{y}{L}\right)-\phi _{K(\lambda )}\left( \frac{y}{L}\right)>0$ (since $\phi _{K(\lambda )}$ is even and negative), from which we also get the smoothness of $L(y,z)$ in $\mathcal{C}_\lambda \cap\Omega$ because $\phi _{K(\lambda )}$ is smooth. The continuity of $L$ up to the boundary follows from the definition and the continuity of $\phi _{K(\lambda )}$.
\end{Proofc}
Using the diffeomorphism from Lemma \ref{Lem_subsol_2} we can define $q=q_\lambda (y,z)\in C^\infty (\mathcal{C}_\lambda \cap \Omega ,\mathbb{R}^2)\cap C(\overline{\mathcal{C}_\lambda\cap \Omega}\setminus \{(0,0)\} ,\mathbb{R}^2)$ as follows
\begin{equation}\label{eq:vector_field}
q(y,z):=
\frac{(-\phi ' _{L(y,z)}(y),1)}{\sqrt{1+|\phi ' _{L(y,z)}(y)|^2}} , \quad (y,z)\in \overline{\mathcal{C}_\lambda \cap \Omega}\setminus \{(0,0)\} , 
\end{equation}
where $\phi ' _{L(y,z)}(y)=\phi ' _{K(\lambda )}\left( \frac{y}{L(y,z)}\right)$. Note that the boundary values of $q$ make sense because of the boundary values of $\phi '_{K(\lambda )}$ by Lemma \eqref{eq:expl_implicite_solution}.
We have the following Lemma
\begin{lem}{\bf(An equation for $q$)\\}\label{Lem:sub_sol_3}
Let $\lambda >1$, $q$ as in \eqref{eq:vector_field} then
\begin{equation}\label{eq:sol_vecttor_field}
-\mathrm{div}(|z|q(y,z))=\lambda , \quad \text{for } (y,z)\in (\mathcal{C}_\lambda \cap \Omega ) .
\end{equation}
\end{lem}
\begin{Proofc}{\bf Proof of Lemma \ref{Lem:sub_sol_3}\\}
All the equations in this proof hold for $(y,z)\in \mathcal{C}_\lambda \cap\Omega$.
Having in mind the diffeomorphism $(y,z)\mapsto (\bar{y},L(y,z))$, with $\bar{y}(y)=y$ from Lemma \ref{Lem_subsol_2}, we can write $q=q(\bar{y},L(y,z))=(q_1(\bar{y},L(y,z)),q_2(\bar{y},L(y,z)))$. Since $|z|=-z$ in $\Omega$, we have
\begin{equation}\label{eq:Lem_subsol__1}
\mathrm{div}_{(y,z)}(|z|q)=|z|\mathrm{div}_{(y,z)}(q)-q_2
\end{equation}
and 
\begin{align}\label{eq:Lem_subsol__2}
\partial _yq_1 & =\partial _{\bar{y}}q_1+\partial _L(q_1)\partial _yL\\
\partial _z q_2 &= \partial _L(q_2) \partial _z L \nonumber .
\end{align}
In order to simplify the notation we set $\psi =\phi _{K(\lambda )}$, then using \eqref{eq:implicit_def_K} we can write $q_1=\frac{-\psi '(y/L)}{\sqrt{1+|\psi '(y/L)|^2}}$ and $q_2=\frac{1}{\sqrt{1+|\psi '(y/L)|^2}}$, from which we can calculate
\begin{align}\label{eq:Lem_subsol__3}
\partial _L q_1 &= \frac{\psi ''\left(\frac{y}{L}\right) \left( \frac{y}{L^2}\right) }{\left(1+|\psi '\left( \frac{y}{L} \right)|^2 \right) ^{3/2}},\\
\partial _L q_2 &=\frac{\psi ''\left(\frac{y}{L}\right)\psi '\left( \frac{y}{L}\right) \left( \frac{y}{L^2}\right)}{\left(1+|\psi '\left( \frac{y}{L}\right)|^2 \right) ^{3/2}}.\nonumber
\end{align}
Differentiating \eqref{eq:implicit_def_K} in $y$ and $z$ we get
\begin{align}\label{eq:Lem_subsol__4}
\partial _z L &=\frac{-1}{\frac{y}{L}\psi '\left( \frac{y}{L}\right)-\psi \left( \frac{y}{L}\right)} ,\\
\partial _y L &= \frac{\psi '\left( \frac{y}{L}\right)}{\frac{y}{L}\psi '\left( \frac{y}{L}\right)-\psi \left( \frac{y}{L}\right)}. \nonumber
\end{align}
Using \eqref{eq:Lem_subsol__3}, \eqref{eq:Lem_subsol__4} we get $\partial _L(q_1)\partial _y L+\partial _L(q_2)\partial _z L=0$ and hence we get from \eqref{eq:Lem_subsol__2} 
\begin{equation}\label{eq:Lem_subsol__5}
\mathrm{div}_{(y,z)}q=\partial _{\bar{y}}q_1=\frac{d}{d\bar{y}}\frac{-\phi _L ' (\bar{y})}{\sqrt{1+|\phi _L '(\bar{y})|^2}}=\frac{-\phi _L ''(\bar{y})}{(1+|\phi _L '(\bar{y})|^2)^{3/2}} .
\end{equation}
Using the fact that $\bar{y}=y$, $z<0$, \eqref{eq:implicit_def_K} and \eqref{eq:Lem_subsol__5}, equation \eqref{eq:Lem_subsol__1} becomes 
\[
-\mathrm{div}_{(y,z)}(|z|q)=-\frac{\phi _L(y)\phi _L ''(y)}{(1+|\phi _L '(y)|^2)^{3/2}}+\frac{1}{\sqrt{1+|\phi _L '(y)|^2}} ,\quad 
\]
and finally using the above equation together with \eqref{eq:eq_of_phi_K_max} and the definition of $\phi _L$ we conclude
\[
-\mathrm{div}_{(y,z)}(|z|q)=\lambda ,\quad (y,z)\in \mathcal{C}_\lambda \cap\Omega .
\]
\end{Proofc}
Note also that by \eqref{eq:Lem_subsol__4} and the boundary conditions of $\phi '_{K(\lambda )}$ we can extend $L\in C^1(\overline{\mathcal{C}_\lambda \cap \Omega}\setminus \{(0,0)\})$.
\begin{lem}{\bf (Bound on the Laplacian)\\}\label{Lem:Lapl_bound}
Let $L$ be as in \eqref{eq:implicit_def_K}, then there are positive constants $C_1,C_2$ such that if
\begin{equation}\label{eq:bound_1_deriv}
C=C(\lambda ):= 1+\left(  \frac{\lambda -1}{K(\lambda )} \right) ^2 
\end{equation}
we have
\begin{align}
 (\partial _yL(y,z)) ^2+(\partial _z L(y,z))^2\leq C\quad  & \text{ for } (y,z)\in \mathcal{C}_\lambda \cap \Omega ,\label{eq:lapl_bound_1_deriv_2} \\
  L(y,z)\Delta _{(y,z)} L(y,z)\leq C_1+C_2\quad & \text{ for } (y,z)\in \mathcal{C}_\lambda \cap \Omega \label{eq:lapl_bound_2_deriv}, 
\end{align}
in particular we have
\begin{equation}\label{eq:bound_lapl}
0\leq \Delta _{(y,z)}(L(y,z))^2\leq 2(C+C_1+C_2) ,\quad (y,z)\in \mathcal{C}_\lambda \cap \Omega .
\end{equation}
\end{lem}
\begin{Proofc}{\bf Proof of Lemma \ref{Lem:Lapl_bound}\\ }
As in the proof of Lemma \ref{Lem:sub_sol_3} we simplify the notation by setting $\psi =\phi _{K(\lambda )}$.\\
\textbf{Step 1. Bound on $\partial _z L$ and $\partial _y L$}\par
By \eqref{eq:second_der} we have $\psi ''>0$ in $(-1,1)$, then, using also \eqref{eq19} we can estimate by the maximum
\begin{equation}\label{eq:lapl_bound_0.1}
\frac{1}{y\psi '(y)-\psi (y)}\leq \frac{\lambda -1}{-K(\lambda )} ,\quad \text{ for } y\in (-1,1).
\end{equation}
Let $(y,z)\in \mathcal{C}_\lambda \cap \Omega$, by the diffeomorphism in Lemma \ref{Lem_subsol_2} we have $\frac{|y|}{L(y,z)}\leq 1$, hence using \eqref{eq:lapl_bound_0.1} and the formula of $\partial _z L$ by \eqref{eq:Lem_subsol__4} we have $|\partial _zL(y,z)|\leq \frac{\lambda -1}{-K(\lambda )}$.

Similarly for $\partial _y L$ given by the formula \eqref{eq:Lem_subsol__4}, since 
\[
\frac{d}{dy}\left( \frac{\psi '(y)}{y\psi '(y)-\psi (y)}\right) =\frac{-\psi (y)\psi ''(y) }{(y\psi '(y)-\psi )^2}>0\, \quad \text{ for } y\in (-1,1),
\]
and $\displaystyle\lim_{y\rightarrow 1}\frac{\psi '(y)}{y\psi '(y)-\psi (y)}=1$, we have $|\partial _yL(y,z)|\leq 1$ for $(y,z)\in \mathcal{C}_\lambda \cap \Omega$. Combining the bounds of $\partial _zL$ and $\partial _yL$ we get \eqref{eq:lapl_bound_1_deriv_2}.\\
\textbf{Step 2. Bound on second derivatives}\par
If we differentiate \eqref{eq:implicit_def_K} twice in $z$ and $y$ respectively and use \eqref{eq:Lem_subsol__4} we get
\begin{equation}\label{eq:lapl_bound_3}
L\partial ^2 _{zz} L=\frac{\psi ''\left( \frac{y}{L}\right) \left( \frac{y}{L}\right) ^2}{\left( \frac{y}{L}\psi '\left( \frac{y}{L}\right) -\psi \left( \frac{y}{L}\right) \right) ^3}
\end{equation}
and
\begin{equation}\label{eq:lapl_bound_4}
L\partial ^2 _{yy} L=\frac{\psi ''\left( \frac{y}{L}\right) \psi ^2\left( \frac{y}{L}\right) }{\left( \frac{y}{L}\psi '\left( \frac{y}{L}\right) -\psi \left( \frac{y}{L}\right) \right)^3} .
\end{equation}
We estimate in $\mathcal{C}_\lambda \cap \Omega $ 
\[
L\partial ^2 _{zz}L\leq \max\left\lbrace \max_{\frac{|y|}{L}\leq \frac{1}{2}} L\partial ^2 _{zz}L, \sup_{\frac{1}{2}<\frac{|y|}{L} <1} L\partial ^2 _{zz}L\right\rbrace.
\]
Using the fact that $y\psi '(y)\geq 0$ and the maximum of $\psi$ by \eqref{eq19} we estimate
\begin{equation}\label{eq:lapl_bound_4.25}
\displaystyle \max_{\frac{|y|}{L}\leq \frac{1}{2}} L\partial ^2 _{zz}L \leq \frac{1}{4}\left( \frac{\lambda}{-K(\lambda )}\right) ^3\displaystyle\max_{|y|\leq \frac{1}{2}}\psi ''(y).
\end{equation}
For $1/2<|y/L|<1$ it is $\psi '\neq 0$ and we can rewrite \eqref{eq:lapl_bound_3} as 
\begin{equation}\label{eq:lapl_bound_4.5}
L\partial ^2 _{zz} L=\frac{\psi ''\left( \frac{y}{L}\right)}{\left|\psi '\left( \frac{y}{L}\right)\right| ^3}\cdot \frac{\left( \frac{y}{L}\right) ^2}{\left( \left| \frac{y}{L}\right|+\frac{-\psi \left( \frac{y}{L}\right)}{\left| \psi '\left( \frac{y}{L}\right)\right|}\right)^3} ,
\end{equation}
and by equation \eqref{eq:second_der} we calculate in the same interval
\begin{equation}\label{eq:lapl_bound_5}
\frac{\psi ''}{|\psi '|^3}=\frac{\lambda \left( \frac{1}{|\psi '|^2}+1\right) ^{3/2}-\frac{1}{|\psi '|^3}-\frac{1}{|\psi '|}}{-\psi} .
\end{equation}
Substituting \eqref{eq:lapl_bound_5} in \eqref{eq:lapl_bound_4.5} and using properties of $\psi$ and the monotonicity of $\psi '$ we get the bound
\begin{equation}\label{eq:lapl_bound_6}
\displaystyle\sup_{\frac{1}{2}<\frac{|y|}{L}<1 }L\partial ^2 _{zz} L\leq \frac{2\lambda ^2}{-K(\lambda )}\left( \frac{1}{\left|\psi '\left( \frac{1}{2}\right)\right| ^2}+1\right) ^{3/2}  .
\end{equation}
Finally by \eqref{eq:lapl_bound_6} and \eqref{eq:lapl_bound_4.25} we get $\displaystyle\sup_{\mathcal{C}_\lambda \cap\Omega} L\partial ^2 _{zz}L \leq C_1$, with $C_1$ a positive constant.
Similarly one can show that $\displaystyle\sup_{\mathcal{C}_\lambda \cap \Omega } L\partial ^2 _{yy}L \leq C_2$ with 
\[
C_2=\max \left\lbrace\frac{\lambda}{-K(\lambda )}\max_{|y|\leq\frac{1}{2}}\psi ''(y), 8\left( \frac{\lambda}{\lambda -1}\right)^2 (-K(\lambda ))\left( \frac{1}{\left| \psi '\left(\frac{1}{2}\right)\right| ^2}+1\right) ^{3/2}\right\rbrace .
\]
\end{Proofc}
\subsection{A subsolution}\label{SecSubsolution}
\begin{figure}[htbp]
     \centering
     \subfloat[][Curves with normal $d^{\mathrm{ext}} _{\lambda }$]{\includegraphics[scale=0.8]{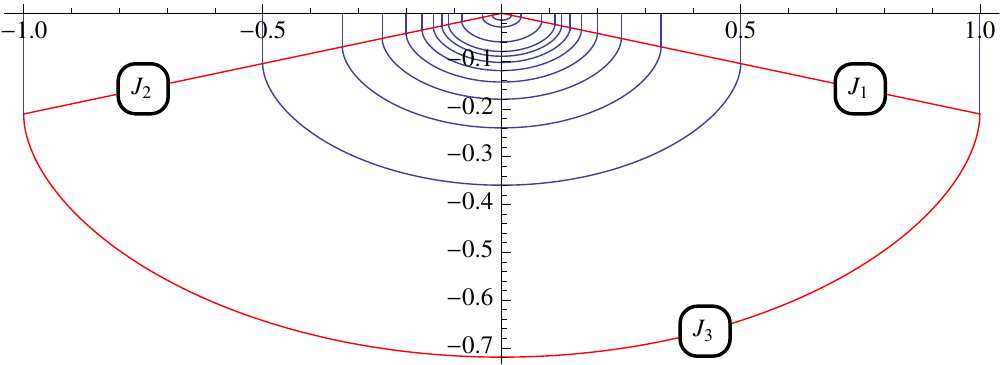}\label{<figure1>}}
     \subfloat[][$u_{\zeta ,\lambda }$]{\includegraphics[scale=0.55]{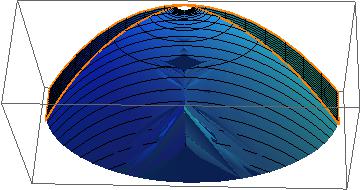}\label{<figure2>}}
     \caption{ }
     \label{fig:avalance}
\end{figure}

\begin{rem}\label{rem:dirac}
Let $\sigma: \Omega _1\cup \Omega _2 \rightarrow\mathbb{R}^2$, with $\Omega _1,\Omega _2\subset \mathbb{R}^2$, two bounded domains with Lipschitz boundary and a common smooth boundary $\partial \Omega $, with surface measure $dS$. Suppose that $\sigma \in \displaystyle\bigcap_{i=1}^{2}(C^1(\Omega _i,\mathbb{R}^2)\cap C(\overline{\Omega _i},\mathbb{R}^2))$, $\mathrm{div} \sigma \in L^2(\Omega _1)\cap L^2(\Omega _2)$, we denote by $\mathrm{Tr}\mid _{\Omega _i} \sigma $, $i=1,2$, the limit value of $\sigma $ from the sides $\Omega _i$ respectively. Then for $\phi \in W^{1,2} _0(\Omega _1 \cup \Omega _2)$ with $\mathrm{supp}(\phi )\cap \partial \Omega \neq \emptyset$ it is 
\begin{equation}\label{eq:dirac_gen}
\int_{\Omega _1\cup\Omega _2} \sigma \cdot \nabla \phi=-\int_{\Omega _1\cup\Omega _2}\mathrm{div}(\sigma )\phi  +\int_{\partial \Omega }n \cdot (\mathrm{Tr}\mid _{\Omega _1} \sigma -\mathrm{Tr}\mid _{\Omega _2} \sigma )\phi\, dS
\end{equation}
where $n$ is the normal to $\partial \Omega$ pointing at the direction of $\Omega _2$.
\end{rem}
We can now construct a subsolution. In what follows we will favour intuition over mathematical elegance, as far as the notation is concerned, and we will instead denote the set $\mathrm{Epi} ^\smallsmile (\lambda )$ defined in \eqref{eq:Epi}, simply by $\{z>\phi _{K(\lambda )}\}$. Let $\zeta >0$, using the diffeomorfism from Lemma \ref{Lem_subsol_2} we can define the continuous function (see Figure \ref{fig:avalance})
	\begin{equation}\label{eq:subsolution_final}
u_{\zeta ,\lambda }(y,z):=\begin{cases} 
-\zeta y^2+\zeta  & \text{ in } \Omega\setminus \mathcal{C}_\lambda ,\\
-\zeta L^2 (y,z) +\zeta &   \text{ in }\mathcal{C}_\lambda  \cap \{z\geq \phi _{K(\lambda )}\} ,\\
0 & \text{ in } \{ z< \phi _{K(\lambda )}\} ,
\end{cases} 
\end{equation}
and for $q_\lambda $ as in \eqref{eq:vector_field} we define
\begin{equation}\label{eq:subsolution_vector_field}
d^{\mathrm{ext}} _{\lambda }(y,z):=\begin{cases} 
\left( -\frac{y}{|y|},0\right) & \text{ in } \Omega\setminus \mathcal{C}_\lambda , \\
q_\lambda (y,z)  &  \text{ in }\mathcal{C}_\lambda  \cap\Omega .
\end{cases} 
\end{equation}
Then we have that $u_{\zeta ,\lambda }\in \mathcal{X}$ with $\partial _z u_{\zeta ,\lambda }(y,0)=0$ for $y\in (-1,1)$ and $d^{\mathrm{ext}} _\lambda\in \Lambda $.
In the set $\mathcal{C}_\lambda  \cap\{ z> \phi _{K(\lambda )}\}$ we have
$\nabla u_{\zeta ,\lambda }=-2\zeta L(\partial _y L,\partial _z L)$,
then using also \eqref{eq:Lem_subsol__4}, \eqref{eq:subsolution_final}, \eqref{eq:subsolution_vector_field}, definition \eqref{eq:vector_field} and the properties of $\phi _{K(\lambda )}$ by Lemma \ref{Lem:explicit_minim} we have that $d^{\mathrm{ext}} _\lambda\cdot \nabla u_{\zeta ,\lambda }=|\nabla u_{\zeta ,\lambda }|$ a.e. in $\Omega$. 
\begin{pro}{\bf (Subsolution)\\}\label{Pro2}
Let $\lambda >1$, then there is $1<\lambda _0<\lambda$ such that for $0<\zeta _0\leq \frac{\lambda -\lambda _0}{2(C+C_1+C_2)}$, with $C,C_1,C_2$ given by Lemma \ref{Lem:Lapl_bound}, the pair $(u_{\zeta _0 ,\lambda _0}, d^{\mathrm{ext}} _{\lambda _0} )$ given by \eqref{eq:subsolution_final}-\eqref{eq:subsolution_vector_field}, is a subsolution of the equation \eqref{eq:sub_differential}.
\end{pro}
\begin{Proofc}{\bf Proof of Proposition \ref{Pro2}\\}
\textbf{Step 1. The subsolution inequalities}\par
We will first show the subsolution inequalities in the set 
\[
\Omega _1\cup \Omega _2 \cup \Omega _3:=\left( \Omega\setminus \overline{\mathcal{C}_{\lambda _0}}\right) \cup\left( \mathcal{C}_{\lambda _0}\cap \{z>\phi _{K(\lambda _0)}\}  \right) \cup \left( \mathcal{C}_{\lambda _0}\setminus \{z\geq\phi _{K(\lambda _0)}\} \right)
\]
where the functions $u_{\zeta _0,\lambda _0},\, d^{\mathrm{ext}} _{\lambda _0}$ are smooth.
Using \eqref{eq:sol_vecttor_field} and \eqref{eq:subsolution_vector_field} we calculate
\begin{equation}\label{eq:Pro2_1}
-\mathrm{div} (|z|d^{\mathrm{ext}} _{\lambda _0}(y,z))=\begin{cases} 
0 & \text{ in } \Omega\setminus \overline{\mathcal{C}_{\lambda _0}} , \\
\lambda _0 &  \text{ in }\mathcal{C}_{\lambda _0}  \cap\Omega .
\end{cases} 
\end{equation}
Also 
\begin{equation}\label{eq:Pro2_2}
-\Delta u_{\zeta _0,\lambda _0}=\begin{cases} 
2\zeta _0 &\text{ in } \Omega\setminus \overline{\mathcal{C}_{\lambda _0}} , \\
0 &  \text{ in }\Omega \cap \{z<\phi _{K(\lambda _0)}\} .
\end{cases} 
\end{equation}
 and using Lemma \ref{Lem:Lapl_bound} we get in $\mathcal{C}_{\lambda _0}\cap \{ z>\phi _{K(\lambda _0)}\}$
\begin{equation}\label{eq:Pro2_3}
-\Delta u_{\zeta _0,\lambda _0} = \zeta _0\Delta L^2 \leq 2\zeta _0 (C+C_1+C_2) 
\end{equation}
If we now combine \eqref{eq:Pro2_1}-\eqref{eq:Pro2_3}, use the fact that the positive constant $2(C+C_1+C_2)$ depends only on $\lambda _0$, we can choose $\zeta _0\leq \frac{\lambda -\lambda _0}{2(C+C_2+C_2)}$($< \lambda /2$ since $C+C_1+C_2>1$ by \eqref{eq:bound_1_deriv}) and get
\begin{equation}\label{eq:Pro2_4}
-\Delta u_{\zeta _0,\lambda _0}-\mathrm{div} (|z|d^{\mathrm{ext}} _{\lambda _0})\leq \lambda \quad \text{ in } \Omega _1 \cup\Omega _2\cup\Omega _3 .
\end{equation}
It remains to show that inequality \eqref{eq:Pro2_4} holds in the rest of $\Omega$. We will use Remark \ref{rem:dirac} for $\sigma =\nabla u_{\zeta _0,\lambda _0}+|z|d^{\mathrm{ext}} _{\lambda _0}$. Note that $\sigma$ is not defined at $(0,0)$ but we still have that it is bounded near $z=0$ by Lemma \ref{Lem:Lapl_bound}.\\
\textbf{Step 2. The Dirac masses}\par
Note that since $\partial _z u_{\zeta _0,\lambda _0} (y,0)=0$ and therefore $\sigma (y,0) =0$, for $y\in(-1,1)\setminus \{(0,0)\}$, in view of \eqref{eq:dirac_gen}, we do not need to take into account the boundary $\{z=0\}$.
We denote by $J=J_1\cup J_2\cup J_3$ the three parts of the boundary of $\Omega_1 \cup \Omega _2\cup \Omega _3$ as in Figure \eqref{<figure1>}. We will show the subsolution inequalities on $J$. We need to estimate for $(i,j)\in \{(1,2),(2,3)\}$, the terms
\begin{equation}\label{eq:Pro2_4.5}
n_j \cdot \left( \mathrm{Tr}\mid _{\Omega _i}\sigma -\mathrm{Tr}\mid _{\Omega _j} \sigma \right)  ,
\end{equation}
where 
  $n _j$ is the normal of the common boundary pointing in the direction of $\Omega _j$. For $J_1$, the right common boundary of $\Omega _1$ and $\Omega _2$ we have $n _1 =\left( \frac{K(\lambda _0)}{\lambda _0},-1\right) $,  using \eqref{eq:proper_explicit_minim} and \eqref{eq:vector_field} one can see that
that $d^{\mathrm{ext}} _\lambda $ is continuous
in $\Omega$, therefore using \eqref{eq:subsolution_final} we get
\begin{equation}\label{eq:Pro2_final1}
 \mathrm{Tr}\mid _{\Omega _1} \sigma -\mathrm{Tr}\mid _{\Omega _2}\sigma = (-2\zeta _0 y,0)+2\zeta _0 L(\partial _y L,\partial _z L)=0 ,\quad (y,z)\in J_1 ,
\end{equation}
where we used
the fact that $y=L$ on $J_1$ and $(\partial _y L,\partial _z L)=(1,0)$ by the Neumann conditions in \eqref{eq:proper_explicit_minim}.
In a similar way we can write \eqref{eq:Pro2_4.5} on $J_2$ as
\begin{equation}\label{eq:Pro2_final2}
-n  _2 \cdot \left( \mathrm{Tr}\mid _{\Omega _1}\sigma -\mathrm{Tr}\mid _{\Omega _2} \sigma \right) =0.
\end{equation}
where $n _2=\left( \frac{-K(\lambda _0)}{\lambda _0},-1\right)$. On $J_3$ we simplify the notation and set $\psi =\phi _{K(\lambda _0)}$, then \eqref{eq:Pro2_4.5} becomes
\begin{align}\label{eq:Pro2_6}
n _3 \cdot \left( \mathrm{Tr}\mid _{\Omega _2}\sigma -\mathrm{Tr}\mid _{\Omega _3} \sigma \right) &=\left( \frac{\psi '}{\sqrt{1+|\psi '|^2}},\frac{-1}{\sqrt{1+|\psi '|^2}}\right) \cdot \left( -2\zeta _0 L(\partial _yL,\partial _z L)\right) \\
&=-2\zeta _0 \frac{\sqrt{1+|\psi '|^2}}{y\psi '-\psi}\leq 0,\nonumber
\end{align}
where in the last equality we used equations \eqref{eq:Lem_subsol__4} and that $L=1$ on $J_3$. 
We can now conclude from estimates \eqref{eq:Pro2_final1}, \eqref{eq:Pro2_final2} and \eqref{eq:Pro2_6}.
\end{Proofc}
\begin{Proofc}{\bf Proof of Theorem \ref{Pro:main_prop} (lower bound)\\}
If we compare the subsolution $u_{\zeta _0,\lambda _0}$ by Proposition \ref{Pro2} with the solution $u_\lambda$ of \eqref{eq:sub_differential} using Proposition \ref{comp-prin}, we get $0\leq u_{\zeta _0,\lambda _0}\leq u_\lambda $ in $\Omega$ for all $\lambda _0\in (1,\lambda )$, hence by definitions \eqref{eq:subsolution_final} and \eqref{eq:Epi} we get
\begin{equation}\label{eq:inclus_1}
\{ u_{\zeta _0,\lambda _0}>0\}=\{z>\phi _{K(\lambda _0)}\}=\mathrm{Epi} ^\smallsmile (\lambda _0)\subset \{u_\lambda >0\} ,\quad \text{ for all } \lambda _0\in (1,\lambda ).
\end{equation}
We set $\psi (y,\lambda )=\phi  _{K(\lambda )}(y)$ for $(y,\lambda )\in [0,1]\times (1,+\infty )$.
By definition \eqref{eq:expl_implicite_solution} we have that $\psi $ satisfies the equation $F(y,\lambda ,\psi (y,\lambda ))=0$ with
\[F:\{(y,\lambda ,z):\, y\in (0,1),\, \lambda \in(1,+\infty ),\, z\in \left( \frac{K(\lambda )}{\lambda-1}, \frac{K(\lambda )}{\lambda}\right) \}\rightarrow \mathbb{R}\]
given by
\[
F(y,\lambda ,z)=K(\lambda )f_\lambda \left( \frac{z}{K(\lambda )}\right)-K(\lambda )f_\lambda \left( \frac{1}{\lambda -1}\right)-y .
\]
The using the formulas \eqref{eq:inverce_funct}, \eqref{eq:K_max} and \eqref{eq:A_lambda} one can check that $F$ is smooth in the domain of it's definition. Since $f' _\lambda \left( \frac{\psi (y,\lambda )}{K(\lambda )}\right) >0$ for $(y,\lambda )\in (0,1)\times (1,+\infty )$ we have by the implicit function theorem that $\psi \in C^\infty ((0,1)\times (1,+\infty ))$.
Since $\phi _{K(\lambda )}$ is even, we get that for fixed $y\in (-1,0)\cup (0,1)$ the function $\phi _{K(\lambda )} (y)$ is continuous in $\lambda$ in $(1,+\infty )$. By the formulas of $\phi _{K(\lambda )} (\pm 1)$, $\phi _{K(\lambda )} (0)$ by Lemma \eqref{Lem:explicit_minim} and the continuity of the function $K(\lambda )$ we get that $\displaystyle\lim_{\lambda _0 \uparrow \lambda}\phi _{K(\lambda _0)} (y)=\phi _{K(\lambda )} (y)$ for all $y\in [-1,1]$. We can now pass to the limit in \eqref{eq:inclus_1} and conclude.
\end{Proofc}

 \subsection{A supersolution}\label{SecSupersolution}
Let $\lambda >1$, $\lambda _1>\lambda $ and $\vartheta ,b,\Pi$ given by \eqref{eq:theta}, \eqref{eq:b}, \eqref{eq:Pi} respectively. Using the diffeomorphism from Lemma \ref{Lem_subsol_2} with $\phi _{K(\lambda _1)}$ in \eqref{eq:implicit_def_K} we can consider sets of the form $\{(y,z)\in \mathcal{C}_\lambda \cap\Omega:\, 1\leq L(y,z)\leq b\}$, where the level set $\{(y,z)\in \mathcal{C}_\lambda \cap\Omega :\,L(y,z)=1\}$ is the graph $\{z=\phi _{K(\lambda _1)}\}$; we will simply denote by $\{1\leq L(y,z)\leq b\}$ these sets.
We define
\begin{equation}\label{eq:supe_1_part}
u_1 ^{\lambda _1}(y,z):=\frac{\lambda _1}{2}(1-y^2),\quad (y,z)\in \Omega ,
\end{equation}
\begin{equation}\label{eq:supersol_2_part}
u_2 ^{\lambda _1,\vartheta}(y,z):= \begin{cases} 
+\infty   &  \text{ in } \{  z>\phi _{K(\lambda _1)} \} , \\
\vartheta(L(y,z)-b)^2 & \text{ in } \{ 1\leq L(y,z)\leq b\} , \\
0 &  \text{ in } \{  b\leq L(y,z)\} .
\end{cases} 
\end{equation}
where we simply write $\vartheta$ for $\vartheta _{\lambda ,\lambda _1}$. Also, we define 
\begin{equation}\label{eq:supersol}
U^{\lambda _1,\vartheta}=\min\{u_1 ^{\lambda _1} ,u_2 ^{\lambda _1,\vartheta }\} ,\quad \text{ in } \Omega .
\end{equation}
\begin{figure}[htbp]
     \centering
     \includegraphics[scale=0.5]{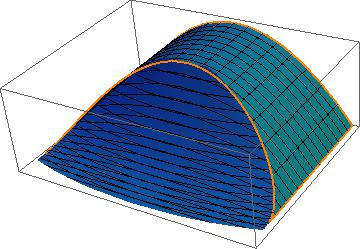}
     \caption{$U^{\lambda _1,\vartheta }$ }
    \label{fig:decreasing_3}
\end{figure}
We note that the intersection of the graphs of the functions $u_2 ^{\lambda _1,\vartheta }$ and $u_1 ^{\lambda _1}$ lies in the domain $\Omega \cap \{L(y,z)<b\}$ and is given by the equation
\begin{equation}\label{eq:intersection}
\vartheta (L(y,z)-b)^2=\frac{\lambda _1}{2}(1-y^2),\quad (y,z)\in \Omega\cap \{L(y,z)<b\} ,
\end{equation}
or else since $L<b$
\[
L(y,z)=b-\sqrt{\frac{\lambda _1}{2\vartheta}(1-y^2)}\geq b-\sqrt{\frac{\lambda _1}{2\vartheta}}=1,
\]
by the definition of $b$. Also, since $\partial _zL<0$ in $\mathcal{C}_\lambda \cap \Omega$ the curve defined by the contour \eqref{eq:intersection} is the graph of a function which lies in fact in the set $\{1\leq L(y,z)\leq b\}$, and therefore, the function $U^{\lambda _1,\vartheta }$ is continuous, see Figure \ref{fig:decreasing_3}.
 For $q_{\lambda _1}$ as in \eqref{eq:vector_field} we define for a.e. $y\in \Omega$ the vector field
\begin{equation}\label{eq:supersol_vector_field}
q^{\mathrm{ext}} _{\lambda _1}(y,z):= \left\{
\begin{aligned}
&\left( -\frac{y}{|y|},0 \right) & \text{ in } &(\{ 1\leq L(y,z)<b-\sqrt{\frac{\lambda _1}{2\vartheta}(1-y^2)} \} \cup\{z>\phi _{K(\lambda _1)}\}) \cap\{ y\neq 0\} \\
&q_{\lambda _1}(y,z) & \text{ in } &\{  b-\sqrt{\frac{\lambda _1}{2\vartheta}(1-y^2)}<L(y,z) \} .
\end{aligned}
\right.
\end{equation}
We have the following Proposition.
\begin{pro}{\bf (Supersolution)\\}\label{pro:super_sol}
Let $\lambda >1$, then the function $U^{\lambda _1,\vartheta}$ defined in \eqref{eq:supersol} is a supersolution of \eqref{eq:sub_differential}. 
\end{pro}
\begin{Proofc}{\bf Proof of Proposition \ref{pro:super_sol}\\}
A straightforward calculation shows that $\nabla U^{\lambda _1,\vartheta}\cdot q^{\mathrm{ext}} _{\lambda _1}=|\nabla U^{\lambda _1,\vartheta} |$, a.e. in $\Omega $. We also have $\partial _z  U^{\lambda _1,\vartheta} (y,0)=\partial _z u_1 ^{\lambda _1} (y,0)=0$.\\
\textbf{Step 1. Supersolution inequalities}\par
It is 
\[
-\Delta U^{\lambda _1,\vartheta}= \left\{
\begin{aligned}
&\lambda _1& \text{ in } &\{1\leq  L(y,z)<b-\sqrt{\frac{\lambda _1}{2\vartheta}(1-y^2)}\}\cup \{z>\phi _{K(\lambda _1)}\} \\
& -2\vartheta ((\partial _y L)^2+(\partial _zL)^2) & &\\
& \quad +2\vartheta   (b-L)(\partial ^2 _{yy}L+\partial ^2 _{zz}L) &\quad \text{ in } &\{  b-\sqrt{\frac{\lambda _1}{2\vartheta}(1-y^2)}<L(y,z)<b \} \\
&0 &\quad \text{ in } &\{ b< L(y,z) \} ,
\end{aligned}
\right.
\]
and as in \eqref{eq:Pro2_1} we have
\[
-\mathrm{div}(|z|q^{\mathrm{ext}} _{\lambda _1} )= \left\{
\begin{aligned}
&0 &\quad \text{ in } &( \{ 1\leq L(y,z)<b-\sqrt{\frac{\lambda _1}{2\vartheta}(1-y^2)} \} \cup\{z>\phi _{K(\lambda _1)}\}) \cap\{ y\neq 0\}\\
& \lambda _1 &\quad \text{ in } &\{  b-\sqrt{\frac{\lambda _1}{2\vartheta}(1-y^2)} <L(y,z)\} .
\end{aligned}
\right.
\]
Therefore if $C$ is as in \eqref{eq:bound_1_deriv}, we have $\vartheta = \frac{\lambda _1-\lambda }{2  C}$ and
\[
-\Delta U_{\lambda _1,\vartheta }-\mathrm{div}(|z|q^{\mathrm{ext}} _{\lambda _1} )\geq \lambda , \quad \text{ in } \Omega \setminus (\{ L(y,z)=b-\sqrt{\frac{\lambda _1}{2\vartheta}(1-y^2)}\} \cup \{0\}\times \left(\frac{K(\lambda _1)}{\lambda _1-1},0 \right) ).
\]
Note that the solution of the equation $L(0,z)=1$ is $z=\frac{K(\lambda _1)}{\lambda _1-1}$. We also note that by \eqref{eq:lapl_bound_3}, \eqref{eq:lapl_bound_4} and Step 2 of the proof of Lemma \ref{Lem:Lapl_bound} we have that $\Delta U_{\theta ,\lambda _1}$ is bounded.\\
\textbf{Step 2. Dirac masses}\par
The discontinuities of the vector fields $\nabla U_{\lambda _1,\vartheta }$ and $q^{\mathrm{ext}} _{\lambda _1} $ lie on the intersection given by the contour \eqref{eq:intersection} and on $\{0\}\times \left(\frac{K(\lambda _1)}{\lambda _1-1},0\right)$. For the second set only the vector field $q^{\mathrm{ext}} _{\lambda _1}$ is discontinuous and the Dirac mass it creates is
\[
|z|(1,0)\cdot ((1,0)-(-1,0))\geq 0.
\]
For the intersection, eq. \eqref{eq:intersection}, we suppress the indices $\lambda _1,\vartheta $ and we write the Dirac mass as
\begin{equation}\label{eq:pro_super_1}
n \cdot \left[ (\nabla u_1-\nabla u_2 )+|z|\left( \frac{\nabla u_1 }{|\nabla u_1 |}-\frac{\nabla u_2 }{|\nabla u_2 |}\right)\right],
\end{equation}
where $n$ is the normal to the intersection pointing at the direction of $\{L(y,z)>b-\sqrt{\frac{\lambda _1}{2\vartheta}(1-y^2)} \}$. Then the $z-$component of $n$ is negative, and since $L_z<0$ by \eqref{eq:Lem_subsol__4} we have
\[
n =\frac{\nabla u_1 -\nabla u_2 }{|\nabla u_1 -\nabla u_2 |} .
\]
Clearly we have $n \cdot (\nabla u_1 -\nabla u_2 )\geq 0$. The second term of \eqref{eq:pro_super_1} is 
\[
\frac{|z|}{|\nabla u_1-\nabla u_2|} \left( |\nabla u_1|  +|\nabla u_2|-\nabla u_1\cdot \nabla u_2\frac{|\nabla u_1|+|\nabla u_2|}{|\nabla u_1||\nabla u_2|} \right)\geq 0
\]
by the Cauchy-Schwartz inequality.
This concludes the proof.
\end{Proofc}
\begin{Proofc}{\bf Proof of Theorem \ref{Pro:main_prop} (upper bound)\\}
We will estimate $\mathrm{supp}\, u$ from above. By Propositions \ref{pro:super_sol} and \ref{comp-prin} we get $0\leq u_\lambda \leq U_{\lambda _1,\vartheta }$ in $\Omega$ and since $\mathrm{supp}\, U_{\lambda _1,\vartheta }=\mathrm{Epi}_\smallsmile (\lambda _1)$ we get the desired estimate.
\end{Proofc}
Let $\lambda ^\star _1=\lambda ^\star _1(\lambda )>\lambda $ be a minimizer of $\Pi (\lambda ,\cdot )$ (see discussion before Theorem \ref{Pro:main_prop}).
In Figure \ref{fig:test2} we give the graph of $\Pi (\lambda ,\lambda ^\star _1)$ for different values of $\lambda$ and in Table \ref{tab:test1} the corresponding minimizers and minimal values. In fact one notices that the difference $\lambda ^\star _1-\lambda$ increases as $\lambda\rightarrow +\infty$, see Figure \ref{fig:increase_diff}.
\begin{center}
\begin{minipage}{.5\textwidth}
  \centering
  \begin{tabular}{lll}
    \toprule
 $\lambda =1.2$ & $\lambda ^\star _1=1.59451$& $\Pi =3.20584$ \\ 
$\lambda =1.4$ & $\lambda ^\star _1=1.84198$& $\Pi =3.66274$ \\ 
$\lambda =1.6$ & $\lambda ^\star _1= 2.09337$&$\Pi =4.16455$ \\
$\lambda =1.8$ & $\lambda ^\star _1=2.34819$& $\Pi =4.69225$\\
    \bottomrule
  \end{tabular}
  \captionof{table}{Optimal $\lambda ^\star _1$}
  \label{tab:test1}
\end{minipage}%
\begin{minipage}{.5\textwidth}
  \centering
  \includegraphics[scale=0.8]{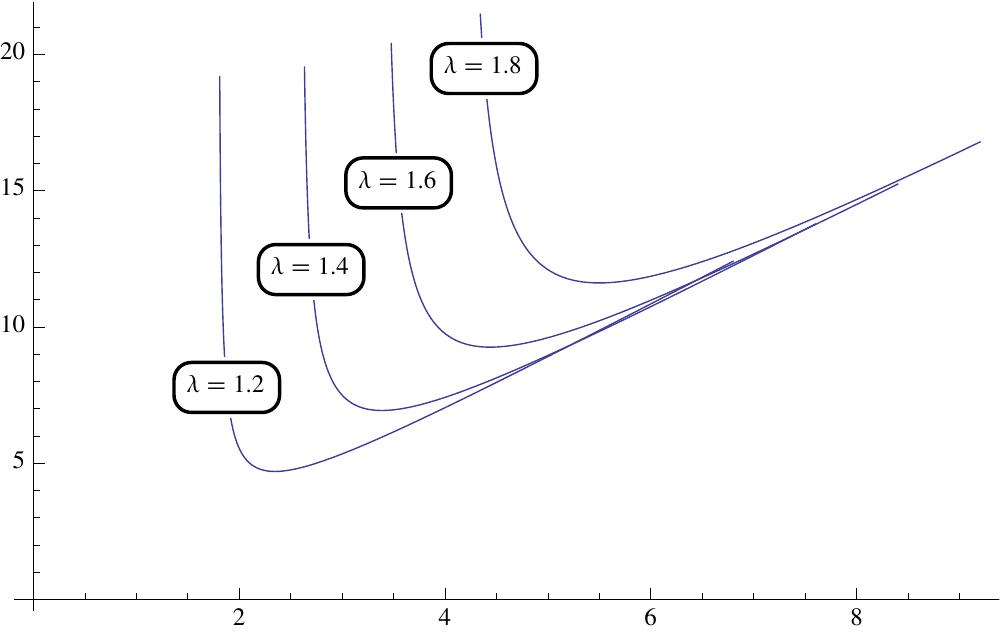}
  \captionof{figure}{$\Pi (\lambda ,\lambda ^\star _1)$ for $\lambda =1.2,\,1.4,\, 1.6,\, 1.8$}
  \label{fig:test2}
\end{minipage}%
\end{center}
\begin{figure}[htbp]
     \centering
     \subfloat[][$\lambda ^\star _1-\lambda$ as $\lambda \rightarrow +\infty$]{\includegraphics[scale=0.8]{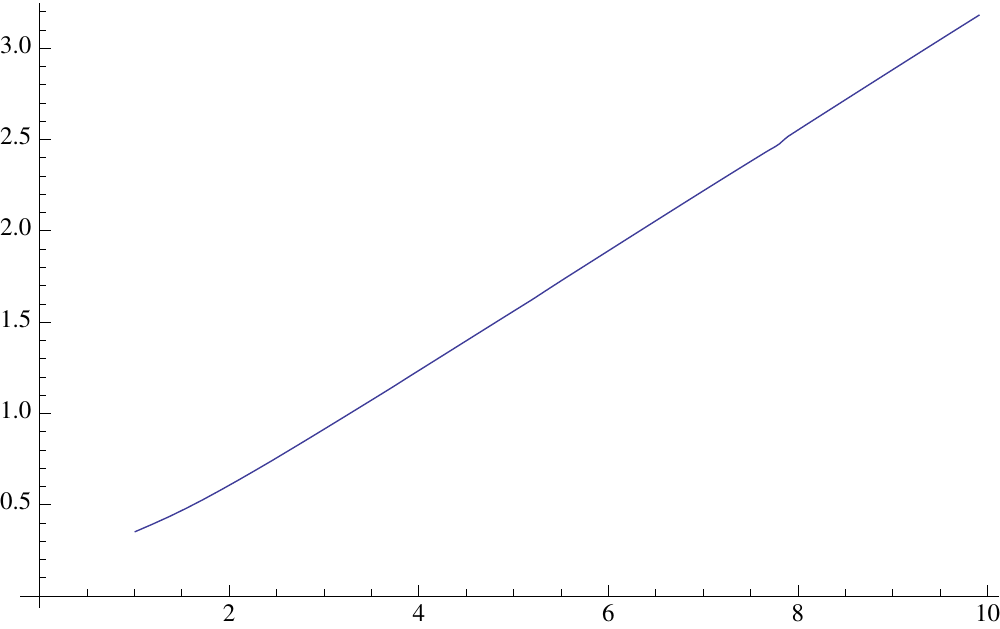}\label{fig:increase_diff}}
     \subfloat[][$\Pi (\lambda ,\lambda ^\star _1)$ as $\lambda \rightarrow +\infty$]{\includegraphics[scale=0.8]{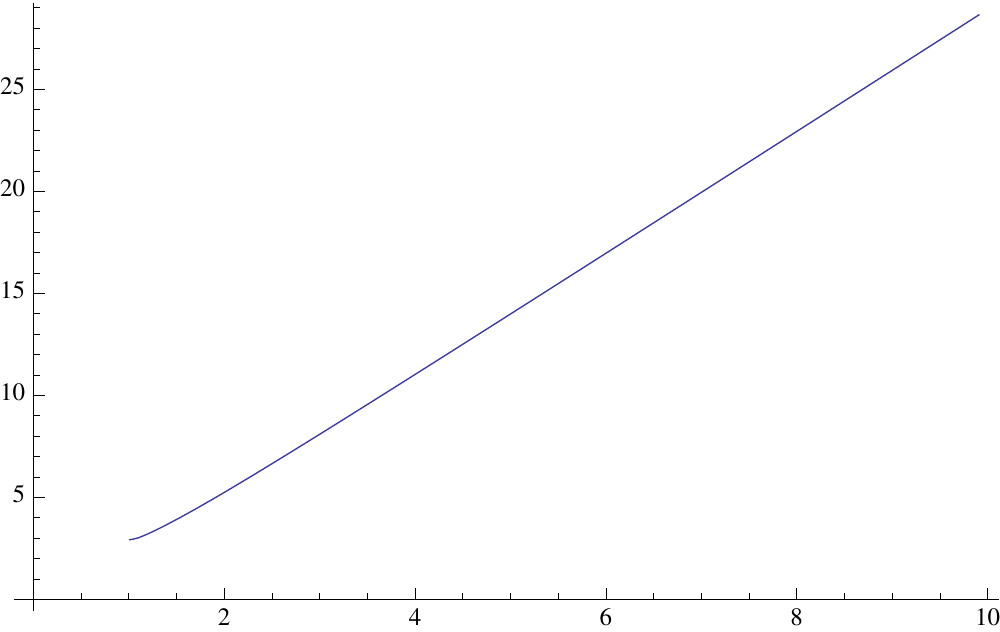}\label{fig:decreasing_2}}
     \caption{ }
     \label{fig:increase_min}
\end{figure}
 \setcounter{section}{0}
\renewcommand{\thesection}{\Alph{section}}
\section{Regularity of $\varepsilon$-minimizers}\label{appendix}
In what follows we will denote by $c$ a generic constant which does not depend on the $\varepsilon$ mentioned in Proposition \ref{Pro:Reg_approx}.
\begin{Proofc}{\bf Proof of Proposition \ref{Pro:Reg_approx} \\}
\textbf{Step 1. Existence/Uniqueness}\par
The uniqueness of the minimizer follows by the strict convexity of the functional or using similar arguments as in the proof of Step 1 of Theorem \ref{theo:main_theo_1} \eqref{itm:E_L_eq_1}. The existence is also similar, in fact the lower semicontinuity of the linear term $-\lambda \int_{\Omega _A}u$ is trivial since the domain $\Omega _A$ is bounded.
We set
\begin{equation}\label{eq:convex_integrad}
F(z,p)=\frac{|p |^2}{2}+|z|\sqrt{\varepsilon ^2 +|p|^2},
\end{equation}
for $(z,p)\in \Omega _A\times \mathbb{R}^2$. It is
\begin{equation}\label{eq:2_der_bd}
|\xi |^2\leq \frac{\partial ^2 F}{\partial p_i\partial p_j}\xi _i\xi _j
\end{equation}
for $\xi =(\xi _1,\xi_2)\in \mathbb{R}^2$, and
\begin{equation}\label{eq:2_der_up_bd}
\left| \frac{\partial ^2 F}{\partial p_i\partial p_j}\right| \leq c\left( 1+\frac{A}{\varepsilon}\right) ,\text{ for all } i,j\in \{1,2\},
\end{equation}
we set $c_2:=c\left( 1+\frac{A}{\varepsilon}\right)$.\\
\textbf{Step 3. Regularity}\par
Since the proof of regularity is standard we are only going to emphasize the particularities of the problem, i.e. the fact that $F$ is only Lipschitz continuous in the $z$ variable. We will simply write $F(z,\nabla u)$ for $F(z,\nabla u(y,z))$. Let $\varphi$ with $\hat{\varphi} \in W^{1,2} _0(\hat{\Omega} _A)$, then equation \eqref{eq:E_L_appr} holds as the first variation of the functional $E^A _{\varepsilon ,\lambda}$. Moreover, using a change of variables one can see that the function $w:=\hat{u} _{\varepsilon ,A}$ satisfies
\begin{equation}\label{eq:E_L_appr_symm}
\int_{\hat{\Omega}_A}\nabla w \cdot \nabla \varphi +|z|\frac{\nabla w \cdot \nabla \varphi}{\sqrt{\varepsilon ^2+|\nabla w |^2}}=\lambda  \int_{\hat{\Omega}_A}\varphi  .
\end{equation}
We study the regularity properties of \eqref{eq:E_L_appr_symm}. Let  $|h|<\mathrm{dist} (\mathrm{supp}\, \varphi ,\partial \hat{\Omega}_A)$, we define $\varphi _{k,h}(y,z):=\varphi ((y,z)-he_k),\, k=1,2$, with $e_k, \, k=1,2$ the unit vectors on the axes $y$ and $z$ respectively. We use $\varphi _{k,h}$ as a test function in \eqref{eq:E_L_appr_symm} and estimate the derivative of the difference quotient   
\begin{equation}\label{eq:diff_quot}
\Delta _h ^k w(y,z) =\frac{w((y,z)+he_k)-w(y,z)}{h} .
\end{equation}
Since the proof is similar we will only present the estimate for $e_2$. Using $\varphi _{2,h}=\varphi _h$ as a test function in \eqref{eq:E_L_appr_symm} and after changing the variables in the integral we get
\begin{equation}\label{eq:Reg_appr_1}
\int_{\hat{\Omega }_A} \partial _{p_i} F(z+h,(\nabla w)_h)\partial _{x_i} \varphi =\lambda \int_{\hat{\Omega}_A} \varphi ,
\end{equation}
where $\partial _{p_i} F=\frac{\partial F}{\partial p_i}$, $(\nabla w)_h(y,z)=\nabla w(y,z+h)$ and $\partial _{x_i}\varphi$, $i=1,2$ is the partial derivative of $\varphi$ in the directions $y,z$ respectively. As usual subtracting \eqref{eq:E_L_appr_symm} from \eqref{eq:Reg_appr_1} we get after a few calculations
\begin{equation}\label{eq:Reg_appr_2}
\int_{\hat{\Omega}_A} \frac{1}{h}(\partial _{p_i} F(z+h,(\nabla w)_h)-\partial _{p_i}F(z+h,\nabla w))\partial _{x_i}\phi =-\int_{\hat{\Omega}_A} \frac{1}{h}(\partial _{p_i} F(z+h,\nabla w)-\partial _{p_i}F(z,\nabla w))\partial _{x_i}\phi . 
\end{equation}
The right hand side of \eqref{eq:Reg_appr_2} can be estimated using the Lipschitz continuity of $\nabla _p F$ in the $z$ variable, we have 
\[
\frac{1}{|h|}|\nabla _p F(z+h,\nabla w)-\nabla _p F(z,\nabla w)|=\frac{|\nabla  w|}{\sqrt{\varepsilon ^2+|\nabla w|^2}}\frac{||z+h|-|z||}{|h|}\leq 1.
\]
It is now a standard process to use \eqref{eq:2_der_bd} and \eqref{eq:2_der_up_bd} in order to bound the quantity $\int |\nabla \Delta _h w|^2$ uniformly in $h$, we have
\begin{equation}\label{eq:Reg_appr_final}
\int_{\Omega '} |\nabla \Delta _h w|^2\leq  2c_3(1+2c_2)\int_{\Omega _A} |\nabla w|^2,
\end{equation}
with $c_3$ a constant independent of $h$ and $\Omega '\subset\subset\Omega ''\subset \subset  \Omega _A$. We then have $w\in W^{2,2}(\Omega '')$ by standard arguments.\\
\textbf{Step 4. Neumann condition}\par
Since $\hat{u} _\varepsilon \in W^{2,2} _{\mathrm{loc}}(\hat{\Omega} _A)$ we can define $\partial _zu_\varepsilon (y,0)$ for a.e. $y\in (-1,1)$ and since $\hat{u}$ is symmetric with respect to $\{z=0\}$, it is in fact $\partial _zu_\varepsilon (y,z)=-\partial _zu_\varepsilon (y,-z)$ for $(y,z)\in \hat{\Omega}_A$; setting $z=0$ we get the desired result.
\end{Proofc}
The constant $c_2$ in the estimate \eqref{eq:Reg_appr_final} depends on $\varepsilon$. Using an argument similar to the proof of \cite[Theorem 3.3.4]{MR1810507} we can show that the second derivative of $u_\varepsilon$ is bounded in $L^2$, uniformly in $\varepsilon$. We have the following Lemma.
\begin{lem}{\bf (Uniform bound on $|\nabla ^2u_\varepsilon |$)\\}\label{lem:unif_bd_2der}
Let $A>0$, $u_\varepsilon $ as in Proposition \ref{Pro:Reg_approx} and $\Omega '\subset \subset \Omega ''\subset\subset \hat{\Omega} _A$. Then there exists a positive constant $C=C(A,\mathrm{dist}(\Omega ',\partial \hat{\Omega} _A))$ such that
\begin{equation}\label{eq:unif_bd_2der}
\int_{\Omega '}|\nabla \partial _{x_i}\hat{u} _\varepsilon |^2\leq C(1+\int_{\hat{\Omega} _A}|\partial _{x_i} \hat{u} _\varepsilon |^2),\quad i=1,2.
\end{equation}
\end{lem}
\begin{Proofc}{\bf Proof of Lemma \ref{lem:unif_bd_2der}\\}
Since the proof is similar to the proof of Proposition \ref{Pro:Reg_approx}, we will only give a sketch of it. We will only show the proof of the estimate \eqref{eq:unif_bd_2der} for $|\nabla \partial _z \hat{u}_\varepsilon |$ because the term with the partial derivative in the $y$ variable is easier to estimate, since the integrand $F$ from \eqref{eq:convex_integrad} does not depend on $y$.
 Let $\varphi$ be a smooth function with compact support in $\Omega ''$; using $\partial _z\varphi$ as a test function in \eqref{eq:E_L_appr_symm} and integrating by parts we can write, using the usual summation convention and the same notation as in the proof of Proposition \ref{Pro:Reg_approx}
 \begin{equation}\label{eq:weak_two_der}
  \int_{\Omega ''}\partial _{p_i}\partial _z F(z,\nabla \hat{u}_\varepsilon )\partial _{x_i}\varphi =0.
 \end{equation}
 Or if we notice that $\partial _z\left( F(z,\nabla \hat{u}_\varepsilon (y,z))\right) =\partial _{\bar{z}}F (\bar{z},\nabla \hat{u}_\varepsilon (y,z))\mid _{\bar{z}=z}+\partial _{p_j}F(z,\nabla \hat{u}_\varepsilon (y,z))\partial _{x_j}\partial _z \hat{u}_\varepsilon (y,z)$ and if $\partial _z|z|=\chi _{(0,+\infty )}-\chi _{(-\infty ,0)}$, we may rewrite \eqref{eq:weak_two_der} as 
  \begin{equation}\label{eq:unif_bd_1}
\int_{\Omega ''} \partial _z|z|\frac{\partial _{x_i} \hat{u}_\varepsilon }{\sqrt{\varepsilon ^2+|\nabla \hat{u}_\varepsilon |^2}}\partial _{x_i}\varphi +\int_{\Omega ''}\left( \partial _{p_i} \partial _{p_j} F(z,\nabla \hat{u}_\varepsilon )\partial _{x_j} \partial _z \hat{u}_\varepsilon \right)\partial _{x_i}\varphi  =0. 
 \end{equation}
As usual we choose a function $\eta\in C^2 _0(\Omega '')$ with $\eta =1$ in $\Omega '$ , $0\leq \eta\leq 1$, $|\nabla \eta |\leq \frac{c}{\mathrm{dist}(\Omega ',\partial \Omega _A)}$ and $||\nabla ^2\eta ||\leq \frac{c}{(\mathrm{dist}(\Omega ',\partial \Omega _A)^2}$. We set $\varphi =\eta ^3\partial _z\hat{u} _\varepsilon $ in \eqref{eq:unif_bd_1}, use the convexity property \eqref{eq:2_der_bd} and the fact that
\[
\partial _{p_i}\partial _{p_j}F(z,\nabla \hat{u}_\varepsilon )\partial _{x_j}\partial _z\hat{u}_\varepsilon =\partial _{x_i}\partial _z\hat{u}_\varepsilon +|z|\partial _z\left( \frac{\partial _{x_i} \hat{u}_\varepsilon }{\sqrt{\varepsilon ^2+|\nabla \hat{u}_\varepsilon |^2}}\right) ,
\]
we get as in the proof of  Proposition \ref{Pro:Reg_approx}
\begin{equation}\label{eq:unif_bd_2}
\int_{\Omega ''}\eta ^3|\nabla \partial _z \hat{u} _\varepsilon |^2 \leq -\int_{\Omega ''}\partial _{x_i}\partial _z \hat{u} _\varepsilon \partial _{x_i}(\eta ^3) \partial _z\hat{u} _\varepsilon -\int_{\Omega ''} \partial _z|z|\partial _{x_i}(\eta ^3) \frac{\partial _{x_i}\hat{u} _\varepsilon }{\sqrt{\varepsilon ^2+|\nabla \hat{u} _\varepsilon |^2}}\partial _z\hat{u}_\varepsilon
\end{equation}
\[
\quad\quad\quad\quad\quad \quad\quad\quad\quad -\int_{\Omega ''}\partial _z |z|\eta ^3 \frac{\partial _{x_i}\hat{u}_\varepsilon }{\sqrt{\varepsilon ^2+|\nabla \hat{u}_\varepsilon |^2}} \partial _{x_i}\partial _z\hat{u}_\varepsilon -\int_{\Omega ''}|z|\partial _z\left( \frac{\partial _{x_i}\hat{u}_\varepsilon}{\sqrt{\varepsilon ^2+|\nabla \hat{u}_\varepsilon|^2}}\right) \partial _{x_i} (\eta ^3) \partial _z\hat{u}_\varepsilon .
\]
The first three terms of the right hand side of \eqref{eq:unif_bd_2} can be estimated as in the proof of  Proposition \ref{Pro:Reg_approx} using Young's inequality, the fact that $\partial _z|z|\leq 1$ and $\frac{|\partial _{x_i}\hat{u}_\varepsilon|}{\sqrt{\varepsilon ^2+|\nabla \hat{u}_\varepsilon |^2}}\leq 1$, for $i=1,2$ uniformly in $\varepsilon$. We will only show the estimate of the last term of \eqref{eq:unif_bd_2}, which we denote by $J$. Integrating by parts $J$ we get
\begin{align}\label{eq:estim_J}
J& =\int_{\Omega ''}\partial _z|z|\frac{\partial _{x_i}\hat{u}_\varepsilon }{\sqrt{\varepsilon ^2+|\nabla \hat{u}_\varepsilon |^2}}\partial _{x_i}(\eta ^3)\partial _z\hat{u}_\varepsilon +\int_{\Omega ''}|z|\frac{\partial _{x_i}\hat{u}_\varepsilon }{\sqrt{\varepsilon ^2+|\nabla \hat{u}_\varepsilon |^2}}\partial _z\partial _{x_i}(\eta ^3)\partial _z\hat{u}_\varepsilon \\
&+ \int_{\Omega ''}|z|\frac{\partial _{x_i}\hat{u}_\varepsilon }{\sqrt{\varepsilon ^2+|\nabla \hat{u}_\varepsilon |^2}}\partial _z(\eta ^3)\partial _z\partial _z\hat{u}_\varepsilon .\nonumber
\end{align} 
It is a standard process now to estimate the right hand side of the above equality using Young's inequality with weight $\gamma >0$, for example the last term of \eqref{eq:estim_J} can be estimated from above by 
\[
c\int_{\Omega ''}\eta ^{1/2} \eta ^{3/2} |\nabla \partial _z\hat{u}_\varepsilon |\leq \tilde{c}(\frac{1}{\gamma} +\gamma \int_{\Omega ''}\eta ^3 |\nabla \partial _z \hat{u}_\varepsilon |^2).
\]
 Finally, putting all the estimates together and choosing $\gamma $ small enough we can absorb the terms $\gamma \int_{\Omega ''}\eta ^3|\nabla \partial _z\hat{u}_\varepsilon |^2$ on the right hand side of \eqref{eq:unif_bd_2} by it's left hand side and by noticing that $\eta =1$ on $\Omega '$ we end up with the desired estimate.
\end{Proofc}
\section*{Acknowledgments}
 The authors would like to thank Professor François Bouchut for the valuable suggestions and the fruitful discussions concerning the connections of the mathematical model with the physical properties. They would also like to thank Professor Alexandre Ern for useful comments on the physical aspects of the problem and Professor Marco Cannone for his helpful suggestions on the presentation of this article.

\def\ocirc#1{\ifmmode\setbox0=\hbox{$#1$}\dimen0=\ht0 \advance\dimen0
  by1pt\rlap{\hbox to\wd0{\hss\raise\dimen0
  \hbox{\hskip.2em$\scriptscriptstyle\circ$}\hss}}#1\else {\accent"17 #1}\fi}

\end{document}